\documentclass[leqno,12pt]{amsart}
\usepackage{amssymb}
\usepackage{amsmath}
\usepackage{enumerate}
\usepackage{amsfonts}
\usepackage{graphicx}

\allowdisplaybreaks
 
\newtheorem{teo}{Theorem}

\newtheorem{lema}{Lemma}
\newtheorem{propo}{Proposition}
\newtheorem{remark}{Remark}

\DeclareMathOperator{\supp}{supp}
\DeclareMathOperator{\pv}{p.v.}

\begin{document}

\title[Harmonic analysis related to Bessel operators]
{Mapping properties of fundamental operators in harmonic analysis related to Bessel operators}

\subjclass[2000]{42C05 (primary), 42C15 (secondary)}
\keywords{Bessel operator, heat-diffusion integral, Poisson integral,
    maximal operator, Riesz transform, square function}

\begin{abstract}
We prove sharp power-weighted $L^p$, weak type and restricted weak
type inequalities for the heat and Poisson integral maximal
operators, Riesz transform and a Littlewood-Paley type square
function, emerging naturally in the harmonic analysis related to
Bessel operators.
\end{abstract}

\author[J. Betancor]{J. J. Betancor}
\address{
    Jorge J{.} Betancor \newline
    Departamento de An\'{a}lisis Matem\'{a}tico\\
    Universidad de la Laguna\\ \newline
    Campus de Anchieta, Avda. Astrof\'{\i}sico Francisco S\'{a}nchez, s/n\\ \newline
    38271 La Laguna (Sta. Cruz de Tenerife), Spain}
\email{jbetanco@ull.es}

\author[E. Harboure ]{E. Harboure}
\address{
    Eleonor Harboure \newline
    IMAL-FIQ\\
    CONICET-Universidad Nacional del Litoral\\ \newline
    Guemes 3450, 3000 Santa Fe, Argentina
} \email{harbour@ceride.gov.ar}

\author[A. Nowak]{A. Nowak}
\address{
      Adam Nowak \newline
      Instytut Matematyki i Informatyki \\
      Politechnika Wroc\l{}awska       \\ \newline
      Wyb{.} Wyspia\'nskiego 27 \\
      50--370 Wroc\l{}aw, Poland
} \email{Adam.Nowak@pwr.wroc.pl}

\author[B. Viviani]{B. Viviani}
\address{
    Beatriz Viviani \newline
    IMAL-FIQ\\
    CONICET-Universidad Nacional del Litoral\\ \newline
    Guemes 3450, 3000 Santa Fe, Argentina
} \email{viviani@ceride.gov.ar}

\thanks{The first-named author was partially supported by MTM2004/05878. 
The second and the fourth-named authors were partially supported by grants from
CONICET and ANPCyT of Argentina and from Universidad Nacional del Litoral.
The third-named author was partially supported by MNiSW Grant N201 054 32/4285.}

\maketitle

\section{Introduction} 

In his monograph \cite{Stein} Stein suggested the study of analogues
of the fundamental operators in the classical harmonic analysis,
such as Riesz transforms, conjugate Poisson integrals, multipliers,
fractional integrals, maximal functions, square functions, in a
context of discrete or continuous expansions with respect to
eigenfunctions of self-adjoint and positive differential operators.
During the last years, this program, or some of its aspects, has
been successfully developed by many authors in various settings.

The study in the framework of Bessel (and also ultraspherical) 
operators was initiated even before \cite{Stein}
by the seminal paper \cite{MS} of Muckenhoupt and Stein. They
introduced the notion of conjugation in the Bessel setting, and their starting 
point was the formulation of suitable Cauchy-Riemann type equations leading to 
a definition of conjugate Poisson integrals. Then the Riesz transform, or
rather the conjugate function mapping according to the terminology used 
in \cite{MS}, emerge as the corresponding boundary value.
After \cite{MS} the Bessel context was investigated by several authors.
In particular, recently Betancor and Stempak
\cite{BS} and Betancor, Buraczewski, Fari\~{n}a, Mart\'{\i}nez and
Torrea \cite{BBFMT,BFMT} obtained some boundedness results for a Riesz
transform and $g$-functions in Bessel settings.

The aim of the present paper is to advance the study of $L^p$
mapping properties of several basic operators related to the harmonic
analysis of Bessel operators. We analyze the
behavior of the maximal operators for the heat and Poisson integrals,
Riesz transform and a $g$-function (see Section \ref{sec:prel} for rigorous definitions
of these objects) associated with the Bessel operator appearing in \cite{MS},
$$
\Delta_\lambda= -\frac{d^2}{dx^2}-\frac{2\lambda}{x} \frac{d}{dx}, \qquad \lambda>-1/2,
$$
which is essentially self-adjoint in
$L^2(\mathbb{R}_+,d\mu_{\lambda})$, with $\mathbb{R}_+=(0,\infty)$ and 
$$
d\mu_{\lambda}(x)=x^{2\lambda}dx, \qquad x>0.
$$
Our main interest is focused on characterizing the power weights
$x^\delta$, for which the abovementioned operators are of strong
type, weak type or restricted weak type $(p,p)$ with respect to the
measure $x^\delta dx$. We shall give a complete description of
such power weights, and in all cases prove the outcomes to be sharp.

Our results are achieved by the nowadays standard method of
splitting the integral kernels into local and global parts, where
\emph{local} is related to a symmetric cone containing the diagonal of
$(0,\infty)\times(0,\infty)$. Following Muckenhoupt and Stein \cite{MS}, we
show that in the local region the operators behave like those derived
from the usual Laplacian, while in the global region they are essentially
controlled by Hardy-type operators. In order to get sharp
results for the range of the power weights, it is necessary to obtain a precise
knowledge of the behavior of the kernels involved. We use several tools in performing
this task, one of them being the local Calder\'on-Zygmund theory established in \cite{NS}.

We point out that in the literature there are some
recent related results regarding the harmonic analysis derived from the Bessel operator
$$
  \widetilde{\Delta}_\lambda = -\frac{d^2}{dx^2} - \frac{\lambda(1-\lambda)}{x^2},	
$$ 
which is essentially self-adjoint in $L^2((0,\infty),dx)$. In particular, 
in \cite{BBFMT,BFMT,BS} Riesz transforms and $g$-functions were studied in this setting.
The results contained in the present paper have counterparts in the framework of
$\widetilde{\Delta}_\lambda$. Moreover, for proving those twin results there is no need to
carry out parallel computations since we may directly take advantage of
the estimates and properties already shown in the
$\Delta_\lambda$ context. Comments sketching how the corresponding results in the
$\widetilde{\Delta}_\lambda$ setting can be concluded will be given along the paper.    
 
Finally, let us give a short account of the previous results concerning the 
operators we investigate. For the Poisson and heat-diffusion integrals, the
unweighted case, with the restriction $\lambda>0$, was studied in \cite{MS} and \cite{BX},
respectively. A $g$-function based on the Poisson kernel was
investigated in \cite{Stem}, where strong type
$(p,p)$ for $p>1$, with respect to the measure $\mu_{\lambda}$, $\lambda>0$, was obtained.
Considering the Riesz transform, in \cite{Anker} a
characterization of the weights for strong type $(p,p)$,
$1<p<\infty$, and weak type $(1,1)$ was given. The approach to this operator was 
analogous to that in \cite{MS}, through conjugate Poisson integrals.
Here we adopt the point of view taken in \cite{BBFMT,BFMT}, and show that for any
$\lambda>-1/2$ the Riesz transform is a principal value integral with a kernel that
satisfies similar estimates to those in \cite{Anker}. Although the
scope of \cite{Anker} for the strong and weak type inequalities is more general than ours, 
we also analyze restricted weak type obtaining new weighted inequalities.
 
The paper is organized as follows. In Section \ref{sec:prel} we 
introduce the main objects of our study and state the main results, Theorems \ref{tmax}-\ref{tg1}. 
There we also gather some general facts and lemmas that will be used throughout the paper. 
The remaining Sections \ref{sec:heat}-\ref{sec:subord}
are separately devoted to the heat integral maximal operator, Riesz transform, $g$-function
and Poisson integral maximal operator, respectively, and the proofs of the main theorems.
In Section \ref{sec:subord} we also take into account a square function related to the Poisson integral.

\section{Preliminaries and statement of results}
\label{sec:prel}

Recall that the standard set of eigenfuctions of the Bessel operator $\Delta_{\lambda}$ consists of
$$
\varphi_z^{\lambda}(x) = (zx)^{-\lambda+1/2}J_{\lambda-1/2}(zx), \qquad x,z>0,
$$
where $J_\nu$ is the Bessel function of the first kind and order $\nu>-1$. 
Indeed, a straightforward computation (cf. \cite[Section 5.2]{Leb}) shows that, for $\lambda>-1/2$,
\begin{equation} \label{eigen}
\Delta_\lambda \varphi^{\lambda}_z = z^2\, \varphi_z^{\lambda}, \qquad z>0.
\end{equation}
Thus the heat kernel associated to $\Delta_\lambda$ is
$$
W_t^\lambda(x,y)=\int_0^\infty e^{-z^2t}
\varphi_z^{\lambda}(x) \varphi_z^{\lambda}(y) \, d\mu_{\lambda}(z), \qquad t,x,y>0.
$$
Computing the last integral (see \cite[p.\,195]{Wat}) leads to
\begin{equation} \label{heat_ker}
W_t^\lambda(x,y)=\frac{(xy)^{-\lambda+1/2}}{2t} e^{-(x^2+y^2)\slash
{4t}} I_{\lambda-1/2} \Big(\frac{xy}{2t}\Big), \qquad t,x,y>0,
\end{equation}
with $I_\nu$ being the modified Bessel function of the first kind and
order $\nu>-1$. Then the heat-diffusion integral of a function $f$ is defined by
$$
W_t^{\lambda}f(x)=\int_0^\infty W_t^\lambda(x,y)f(y) \, d\mu_{\lambda}(y), \qquad t,x>0.
$$
Denote by $W_{*}^{\lambda}$ the corresponding maximal operator,
$$
W_{*}^{\lambda}f = \sup_{t>0} |W_t^{\lambda}f|.
$$

Our result concerning $W_*^{\lambda}$ is the following, see also Figure \ref{fig:max} below.
\begin{teo} \label{tmax}
Let $\lambda > -1\slash 2$, $1\le p < \infty$, $\delta \in \mathbb{R}$.
Then the maximal operator $W_*^{\lambda}$, considered on the measure space
$(\mathbb{R}_+,x^{\delta}dx)$, has the following mapping properties:
\begin{itemize}
\item[(a)] 
$W_*^{\lambda}$ is of strong type $(p,p)$ if and only if
$p>1 \; \textrm{and} \; -1<\delta<(2\lambda+1)p-1;$
\item[(b)]
$W_*^{\lambda}$ is of weak type $(p,p)$ if and only if
$-1<\delta<(2\lambda+1)p-1 \; \textrm{or} \;\;  \delta=2\lambda;$
\item[(c)]
$W_*^{\lambda}$ is of restricted weak type $(p,p)$ if and only if
$-1<\delta\le(2\lambda+1)p-1.$
\end{itemize}
Moreover, $W_*^{\lambda}$ is of strong type $(\infty,\infty)$. 
\end{teo}

\begin{figure}[ht]
\includegraphics[scale=1]{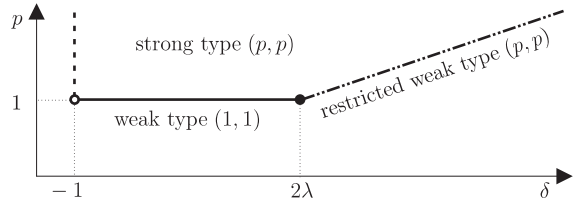}
\caption{Mapping properties of $W_*^{\lambda}$ (fixed $\lambda=1$).}\label{fig:max}
\end{figure}

According to \eqref{eigen}, the Poisson kernel is given by
$$
P_t^\lambda(x,y)=\int_0^\infty e^{-zt}
\varphi_z^{\lambda}(x) \varphi_z^{\lambda}(y)\,d\mu_{\lambda}(z), \qquad t,x,y >0,
$$
and the Poisson integral of a function $f$ is defined as
$$
P_t^{\lambda}f(x)=\int_0^\infty P_t^\lambda(x,y)f(y)\, d\mu_{\lambda}(y), \qquad t,x>0.
$$
Noteworthy, also the Poisson kernel can be computed explicitly, and the resulting
expression contains the Gauss hypergeometric function 
${_2F_1}$ (see Section \ref{sec:subord} for details).

Exactly the same mapping properties as for $W^{\lambda}_*$ turn out to be in force
for the Poisson integral maximal operator $P^{\lambda}_* f = \sup_{t>0} |P^{\lambda}_t f|$.

\begin{teo} \label{SPoi}
Let $\lambda > -1\slash 2$, $1\le p < \infty$, $\delta \in
\mathbb{R}$. Then the maximal operator $P_*^{\lambda}$, considered
on the measure space $(\mathbb{R}_+,x^{\delta}dx)$, has the following mapping properties:
\begin{itemize}
\item[(a)] 
$P_*^{\lambda}$ is of strong type $(p,p)$ if and only
if $ p>1 \; \textrm{and} \; -1<\delta<(2\lambda+1)p-1;$
\item[(b)] 
$P_*^{\lambda}$ is of weak type $(p,p)$ if and only if
$-1<\delta<(2\lambda+1)p-1 \; \textrm{or} \;\;  \delta=2\lambda;$ 
\item[(c)] 
$P_*^{\lambda}$ is of restricted weak type $(p,p)$ if
and only if $ -1<\delta\le(2\lambda+1)p-1.$
\end{itemize}
Moreover, $P_*^{\lambda}$ is of strong type $(\infty,\infty)$. 
\end{teo}

We now pass to the Riesz transform. The Bessel operator can be represented as
$$
\Delta_\lambda=D^*D,
$$
where $D=\frac{d}{dx}$ is the usual derivative and
$D^*=-x^{-2\lambda}\frac{d}{dx}x^{2\lambda}$ is the formal adjoint of $D$
in $L^2(\mathbb{R}_+,d\mu_{\lambda})$. This factorization suggests
the following system of Cauchy-Riemann type equations:
$$
\frac{\partial}{\partial t}P_t^{\lambda}f(x)=-D_x^*Q^{\lambda}_tf(x), \qquad
\frac{\partial}{\partial t}Q_t^{\lambda}f(x)=-D_x P_t^{\lambda}f(x),
$$
with $Q_t^{\lambda}f$ being a suitably defined conjugate Poisson integral,
$$
Q_t^{\lambda}f(x)=\int_0^\infty Q_t^\lambda(x,y)f(y)\, d\mu_{\lambda}(y), \qquad t,x>0.
$$
For $\lambda>0$ the conjugate Poisson kernel that is consistent with the Cauchy-Riemann
type equations has the form
$$
Q_t^\lambda(x,y)=-\frac{2\lambda}{\pi}\int_0^\pi 
\frac{(x-y\cos\theta)(\sin\theta)^{2\lambda-1}}{(x^2+y^2+t^2-2xy\cos\theta)^{\lambda+1}}d\theta,
\qquad t,x,y>0.
$$
Then the Riesz transform $R_\lambda f$ emerges in a natural way as the boundary value of $Q_t^{\lambda}f$,
$$
R_\lambda f(x)=\lim_{t\to 0^+}Q_t^{\lambda}f(x).
$$
This is the classical way of defining $R_{\lambda}$ used by Muckenhoupt and Stein \cite{MS}.
It is known that for each $f\in L^p(\mathbb{R}_+,d\mu_{\lambda})$, $1\le p<\infty$,
the above limit exists for almost every $x>0$.

Nevertheless, our approach to the Riesz transform is more direct, and the definition is
based on a singular integral representation.
By the results of \cite{BBFMT}, both definitions are consistent when $\lambda >0$.

In agreement with a general philosophy, formally the Riesz transform
$R_\lambda$ related to $\Delta_{\lambda}$ has the form
\begin{equation}\label{I1}
R_\lambda f=D\Delta_\lambda^{-1/2}f.
\end{equation}
This becomes rigorous provided that
$f\in C_c^\infty (\mathbb{R}_+)$ and $\lambda>0$, with the potential
operator $\Delta_\lambda^{-1/2}$ expressed in terms of the Poisson integral,
$$
\Delta_\lambda^{-1/2}f(x)=\int_0^\infty P_t^{\lambda}f(x)\,dt, \qquad x>0,
$$
see \cite{BBFMT}. In the present paper, in contrast with \cite{MS,BBFMT,BFMT}, 
we consider the Riesz transform $R_\lambda$ for
the full range $\lambda>-1/2$. Moreover, to define strictly the
operator $\Delta_\lambda^{-1/2}$ we use the heat integral rather than the Poisson one,
$$
\Delta_\lambda^{-1/2}f(x)=\frac{1}{\sqrt{\pi}}\int_0^\infty
\big( W_t^{\lambda}f(x)-\chi_{\{\lambda\le 0\}}W_t^{\lambda}f(0) \big) \frac{dt}{\sqrt{t}}, \qquad x>0.
$$
It will be shown in Section \ref{sec:riesz} that the limit
$W_t^{\lambda}f(0) = \lim_{x\to 0^+} W_t^{\lambda}f(x)$ exists for each $t>0$ and
$\Delta_\lambda^{-1/2}f(x)$ is well defined for $x>0$, provided that $f
\in C^{\infty}_c(\mathbb{R}_+)$. Note that for $-1\slash 2 < \lambda
\le 0$ we have to consider \emph{compensated potentials} in order to
ensure convergence of the defining integral. Then the Riesz
transform $R_\lambda$ of $f\in C_c^\infty(\mathbb{R}_+)$ is defined
by \eqref{I1}. Moreover, for $f\in C_c^\infty(\mathbb{R}_+)$,
$$
R_\lambda f(x)=\pv \int_0^\infty R_\lambda(x,y)f(y)\, d\mu_{\lambda}(y), \qquad x >0,
$$
with the Riesz transform kernel
$$
R_\lambda(x,y)=\frac{1}{\sqrt{\pi}}\int_0^\infty
\frac{\partial}{\partial x}W_t^\lambda(x,y) \frac{dt}{\sqrt{t}}, \qquad x,y >0,\quad x\neq y.
$$
All the details will be given in Section \ref{sec:riesz}. Now we
state the  boundedness properties of $R_\lambda$, see also
Figure \ref{fig:riesz} below. Notice that $R_{\lambda}$ behaves better than the maximal operators.

\begin{teo} \label{triesz}
Let $\lambda > -1\slash 2$, $1\le p < \infty$, $\delta \in \mathbb{R}$.
Then the Riesz transform $R_{\lambda}$, considered on the measure space
$(\mathbb{R}_+,x^{\delta}dx)$, has the following mapping properties:
\begin{itemize}
\item[(a)] 
$R_{\lambda}$ is of strong type $(p,p)$ if and only if
$p>1 \; \textrm{and} \; -1-p < \delta < (2\lambda+1)p-1;$
\item[(b)]
$R_{\lambda}$ is of weak type $(p,p)$ if and only if
$-1-p<\delta<(2\lambda+1)p-1 \; \textrm{or} \;\; \delta \in \{-2, 2\lambda\};$
\item[(c)]
$R_{\lambda}$ is of restricted weak type $(p,p)$ if and only if
$-1-p\le\delta\le(2\lambda+1)p-1.$
\end{itemize}
\end{teo}

\begin{figure}[ht]
\includegraphics[scale=1]{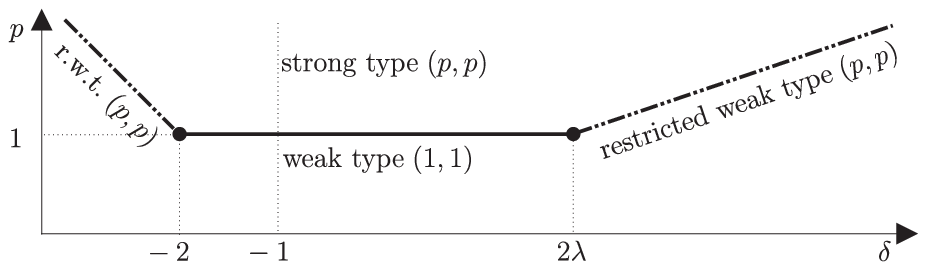}
\caption{Mapping properties of $R_{\lambda}$ (fixed $\lambda=1$).}\label{fig:riesz}
\end{figure}

We now briefly comment on the adjoint Riesz transform $R_{\lambda}^{*}$.
This operator is intimately connected with $R_{\lambda}$, see for instance the identity
\eqref{id_riesz} below. In \cite[(16.8)]{MS} it was shown that for $\lambda>0$
\begin{equation} \label{R1}
R_\lambda f=-xh_{\lambda+1/2}\big(y^{-1}h_{\lambda-1/2}(f)\big),
\qquad f\in L^2(\mathbb{R}_+,d\mu_{\lambda}),
\end{equation}
where $h_\nu$, $\nu>-1$, denotes the \emph{modified Hankel transform},
$$
h_\nu(f)(x)= \int_0^\infty \varphi_x^{\nu+1\slash 2}(y)f(y)\, d\mu_{\nu+1\slash 2}(y)
= \int_0^\infty (xy)^{-\nu}J_\nu(xy)f(y)y^{2\nu+1}dy, \qquad x>0.
$$
It may be proved that (\ref{R1}) holds in fact for
all $\lambda>-1/2$. Since $h_\nu$ is an isometry in $L^2(\mathbb{R}_+,d\mu_{\lambda})$ 
and the Parseval type identity
$$
\int_0^\infty h_\nu(f)(x)g(x)x^{2\nu+1}dx= \int_0^\infty
f(x)h_\nu(g)(x)x^{2\nu+1}dx, \qquad f,g\in
L^2(\mathbb{R}_+,d\mu_{\lambda}),
$$
holds for every $\nu>-1$ (cf. \cite{BS2}), the adjoint of $R_\lambda$ is given by
$$
R_\lambda^* f=-h_{\lambda-1/2}\big(yh_{\lambda+1/2}(x^{-1}f)\big),
\qquad f\in L^2(\mathbb{R}_+,d\mu_{\lambda}).
$$
Further, since $h_\nu^{-1}=h_\nu$ in
$L^2(\mathbb{R}_+,d\mu_{\lambda})$ for all $\nu>-1$, it becomes clear that
\begin{equation} \label{id_riesz}
R_\lambda^* R_\lambda f=R_\lambda R_\lambda^* f=f,
\qquad f \in L^2(\mathbb{R}_+,d\mu_{\lambda}), \quad \lambda > -1\slash 2.
\end{equation}
On the other hand, it is remarkable that $R_\lambda^*$ emerges as
the Riesz transform naturally associated with the Bessel type operator
$$
DD^* = \Delta_{\lambda} + \frac{2\lambda}{x^2}
    = -\frac{d^2}{dx^2} - \frac{2\lambda}{x} \frac{d}{dx} + \frac{2\lambda}{x^2}.
$$
Namely, formally we have $R_{\lambda}^* = D^* (DD^*)^{-1\slash 2}$.
This identity can be given a strict meaning, but we shall not go
into details here to avoid confusion with the line of thought of the paper.

The results concerning $R_{\lambda}^{*}$ can be summarized as follows.

\begin{propo} \label{triad}
Let $\lambda > -1\slash 2$, $1\le p < \infty$, $\delta \in \mathbb{R}$.
Then the adjoint Riesz transform $R_{\lambda}^*$, considered on the measure space
$(\mathbb{R}_+,x^{\delta}dx)$, has the following mapping properties:
\begin{itemize}
\item[(a)] 
$R_{\lambda}^*$ is of strong type $(p,p)$ if and only if
$p>1 \; \textrm{and} \; -1<\delta<2(\lambda+1)p-1;$
\item[(b)]
$R_{\lambda}^*$ is of weak type $(p,p)$ if and only if
$-1<\delta<2(\lambda+1)p-1 \; \textrm{or} \;\;  \delta=2\lambda+1;$
\item[(c)]
$R_{\lambda}^*$ is of restricted weak type $(p,p)$ if and only if
$-1<\delta\le 2(\lambda+1)p-1.$
\end{itemize}
Moreover,
$$
R_\lambda^* R_\lambda f=R_\lambda R_\lambda^* f=f, \qquad f\in L^p(\mathbb{R}_+,x^{\delta}dx),
$$
provided that $p>1$ and $-1<\delta<(2\lambda+1)p-1$.
\end{propo}

Observe that $(a)$ above can be directly deduced from the strong
type result for the Riesz transform. For the other items we may
apply the same arguments as for $R_{\lambda}$ because, as we shall see, they rely on
pointwise estimates of the kernel $R_{\lambda}(x,y)$.
Since the kernel of $R_{\lambda}^{*}$ is $R_{\lambda}(y,x)$, we
easily obtain the same kind of estimates for the adjoint Riesz transform.

Finally, consider the Littlewood-Paley type square function 
$$
g_{\lambda}(f)(x) =\bigg(\int_0^\infty
t\Big|\frac{\partial}{\partial t} W_t^{\lambda}f(x)\Big|^2dt\bigg)^{1/2}.
$$
We prove the following boundedness properties of $g_{\lambda}$.

\begin{teo} \label{tg1}
Let $\lambda > -1\slash 2$, $1\le p < \infty$, $\delta \in \mathbb{R}$. 
Then the square function $g_{\lambda}$, considered on
the measure space $(\mathbb{R}_+,x^{\delta}dx)$, has the following mapping properties:
\begin{itemize}
\item[(a)] 
$g_{\lambda}$ is of strong type $(p,p)$ if and only if
$p>1 \; \textrm{and} \; -1<\delta<(2\lambda+1)p-1;$
\item[(b)]
$g_{\lambda}$ is of weak type $(p,p)$ if and only if $
-1<\delta<(2\lambda+1)p-1 \; \textrm{or} \;\;  \delta=2\lambda;$
\item[(c)]
$g_{\lambda}$ is of restricted weak type $(p,p)$ if and only if 
$-1<\delta\le(2\lambda+1)p-1.$
\end{itemize}
\end{teo}

Notice that the behavior of $g_{\lambda}$ is exactly the same as that of the
maximal operators $W^{\lambda}_{*}$ and $P_{*}^{\lambda}$, see Figure \ref{fig:max}.
Note also that the Riesz transform and its adjoint, taken into account individually, 
behave better than the maximal operators. 
However, if considered simultaneously, they lead back to precisely the same mapping properties
as those of $W^{\lambda}_{*}$ and $P_{*}^{\lambda}$.
In addition, it is worth to mention that the maximal operators are
bounded on $L^p(\mathbb{R}_+,x^{\delta} d\mu_{\lambda})$ for given $1<p<\infty$ and 
$\lambda>-1\slash 2$ if and only if $x^{\delta} \in A_p^{\lambda}$; here
$A_p^{\lambda} = A_p(\mathbb{R}_+,d\mu_{\lambda})$
denotes the Muckenhoupt class of $A_p$ weights associated
with the homogeneous space $(\mathbb{R}_+,|\cdot|,d\mu_{\lambda})$.

As it was already indicated in the Introduction, several harmonic analysis operators 
associated with the Bessel operator
$$
\widetilde{\Delta}_\lambda= x^{\lambda} \Delta_{\lambda} x^{-\lambda} =
-x^{-\lambda}Dx^{2\lambda}Dx^{-\lambda} = \widetilde{D}^* \widetilde{D},
\qquad \widetilde{D} = x^{\lambda} D x^{-\lambda},
$$
were studied earlier, usually with the assumption $\lambda>0$. We now explain
how our present results are related to those in \cite{BBFMT,BFMT}.
Recall that $\widetilde{\Delta}_\lambda$ is associated with
the Lebesgue measure space $(\mathbb{R}_+,dx)$. Consider the
multiplication operator $V_{\lambda}f(x) = x^{-\lambda}f(x)$, which
is obviously an isometry between $L^2(\mathbb{R}_+,dx)$ and
$L^2(\mathbb{R}_+,d\mu_{\lambda})$. The essential observation is
that $V_{\lambda}$ intertwines all relevant operators in both
settings in question and also establishes an $L^2$ equivalence
between them. More precisely, distinguishing with tildes
appropriately defined objects from the $\widetilde{\Delta}_{\lambda}$ setting, we have
$\widetilde{\Delta}_{\lambda} = V_{-\lambda} \Delta_{\lambda}
V_{\lambda}$, $\widetilde{D} = V_{-\lambda} D V_{\lambda}$,
$\widetilde{W}_*^\lambda= V_{-\lambda}W_*^\lambda V_\lambda$,
$\widetilde{P}_*^\lambda= V_{-\lambda}P_*^\lambda V_\lambda$,
$\widetilde{R}_{\lambda} = V_{-\lambda} R_{\lambda} V_{\lambda}$,
$\widetilde{g}_{\lambda} = V_{-\lambda} g_{\lambda} V_{\lambda}$.
Consequently, we can deduce from Theorems \ref{tmax}-\ref{tg1}
the strong boundedness results in $L^p(\mathbb{R}_+,x^{\delta}dx)$
for the corresponding operators in the $\widetilde{\Delta}_\lambda$
context. Furthermore, by applying the same procedures as those
employed in the proofs of the theorems just mentioned, one can also
obtain the desired weak type and restricted weak type mapping
properties, as well as their sharpness. These results, appropriately
stated and justified, complement those from \cite{BBFMT,BFMT}. After
proving Theorems \ref{tmax}, \ref{triesz} and \ref{tg1} we provide
remarks concerning the boundedness properties of the
operators associated with $\widetilde{\Delta}_\lambda$.

An important ingredient of proofs contained in the following sections
are basic differential and asymptotic properties of the modified Bessel
function $I_\nu$. Those listed below can be found for instance in \cite{Wat} or \cite{Leb}.

One of possible definitions of $I_{\nu}(z)$ for, say, $\nu>-1$ and $z>0$ is
\begin{equation} \label{P1}
I_\nu(z)=\sum_{n=0}^\infty\frac{(z/2)^{2n+\nu}}{\Gamma(n+1)\Gamma(n+\nu+1)}.
\end{equation}
A straightforward analysis of the above series shows that
\begin{equation} \label{P2}
\frac{d}{dz}\big(z^{-\nu}I_\nu(z)\big)=z^{-\nu}I_{\nu+1}(z), \qquad z>0, \quad \nu>-1,
\end{equation}
and
\begin{equation} \label{P3}
\lim_{z\to 0^+}z^{-\nu}I_\nu(z)=\frac{1}{2^\nu\Gamma(\nu+1)}, \qquad \nu>-1.
\end{equation}
Furthermore, for $z>0$ and $\nu>-1$, we have the following asymptotic representation:
\begin{equation} \label{P4}
I_\nu(z)=\frac{e^z}{\sqrt{2\pi z}}\bigg(\sum_{k=0}^n
(-1)^k [\nu,k](2z)^{-k}+\mathcal{O}\big(z^{-n-1}\big)\bigg), \qquad n = 0,1,2,\ldots,
\end{equation}
the coefficients being given by $[\nu,0]=1$ and
$$
[\nu,k]=\frac{(4\nu^2-1)(4\nu^2-3^2)\cdot\ldots\cdot(4\nu^2-(2k-1)^2)}{2^{2k}\Gamma(k+1)},
\qquad k=1,2,\ldots.
$$

Objects that will frequently appear in our estimates are the Hardy type operators
\begin{align*}
H_0^\eta f(x) & = x^{-\eta-1}\int_0^x f(y)y^\eta \, dy, \qquad x>0,\\
H_\infty^\eta f(x) & = x^\eta\int_x^\infty f(y)y^{-\eta-1}\, dy, \qquad x>0,
\end{align*}
considered for $\eta > -1$. The relevant mapping properties of $H_0^{\eta}$ and $H_{\infty}^{\eta}$
are gathered in \cite[Lemmas 3.1 and 3.2]{CH} or 
\cite[Lemmas 3 and 4]{HSTV}, see also references given there.
For the sake of completeness and reader's convenience, we reproduce them below.

\begin{lema} \label{Hardy0} 
Let $\eta>-1$ and consider $H_0^\eta$ on the measure space $(\mathbb{R}_{+},x^{\delta}dx)$. Then
\begin{itemize}
\item[(a)] 
$H_0^\eta$ is of strong type $(p,p)$ when $1<p\le\infty$ and $\delta<p(\eta +1) -1$;
\item[(b)]
$H_0^\eta$ is of weak type $(1,1)$ if $\delta\le \eta$;
\item[(c)]
$H_0^\eta$ is of restricted weak type $(p,p)$ when $1<p<\infty$ and $\delta = p(\eta +1) -1$.
\end{itemize}
\end{lema}

\begin{lema} \label{Hardyinfty} 
Let $\eta>-1$ and consider $H_{\infty}^\eta$ on the measure space $(\mathbb{R}_{+},x^{\delta}dx)$. Then
\begin{itemize}
\item[(a)]
$H_\infty^\eta$ is of strong type $(p,p)$ when $1<p<\infty$ and $-\eta p-1<\delta$;
\item[(b)]
$H_\infty^\eta$ is of strong type $(\infty,\infty)$ for any $\delta\in \mathbb{R}$
if only $\eta > 0$;
\item[(c)]
$H_\infty^\eta$ is of weak type $(1,1)$ when $-\eta-1\le \delta$ ($<$ if $\eta=0$);
\item[(d)] 
$H_\infty^\eta$, $\eta \neq 0$, is of restricted weak type $(p,p)$ when $1<p<\infty$ and 
$\delta = -\eta p -1$.
\end{itemize}
\end{lema}

Another object that will be used throughout is the Gauss-Weierstrass kernel
$$
\mathcal{W}_t(x,y)=\frac{1}{\sqrt{4\pi t}} e^{-(x-y)^2\slash 4t}, \qquad t>0, \quad x,y \in \mathbb{R}.
$$
The associated heat integral,
$$
\mathcal{W}_tf(x) = \int_{\mathbb{R}} \mathcal{W}_t(x,y) f(y) \, dy, \qquad x \in \mathbb{R},
$$
represents the classical heat semigroup $\{\mathcal{W}_t\}_{t>0}$ on the real line.
A crucial argument used repeatedly below relies on a comparison of various operators from the
Bessel setting with the corresponding well-known operators related to $\mathcal{W}_t$.

To establish weighted $L^p$ mapping properties of certain auxiliary operators
appearing in the proofs of Theorems \ref{triesz} and
\ref{tg1}, we shall use the following result that can be proved by applying the
local version of the Calder\'on-Zygmund operator theory on the real line 
(or rather its vector-valued variant), that was developed
recently by Nowak and Stempak \cite[Section 4]{NS}. 
It is remarkable that the results from \cite{NS} remain valid in a vector-valued setting.

\begin{lema} \label{CZT}
Let $(B,\|\cdot\|)$ be a separable Banach space.
Assume that $T$ is a local vector-valued Calder\'on-Zygmund operator, i.e. $T$
is a bounded operator from $L^2(\mathbb{R_+},dx)$ into the Lebesgue-Bochner space
$L^2_B(\mathbb{R}_+,dx)$ such that
$$
Tf(x) = \int_{x\slash 2}^{2x} K(x,y)f(y) \,dy, \qquad \textrm{a.e.} \;\; x\notin \supp f,
	\quad f \in C_c^{\infty}(\mathbb{R}_+),
$$ 
where the $B$-valued kernel is weakly measurable and satisfies the standard estimates
$$
\|K(x,y)\| \le \frac{C}{|x-y|}, \qquad \|\nabla_{\! x,y} K(x,y) \| \le \frac{C}{|x-y|^2},
$$
in the 'local' region $0<x\slash 2 < y < 2x$, $x\neq y$.

Then, for each $\lambda \in \mathbb{R}$, the operator $S_{\lambda}$ defined by
$$
S_{\lambda}f(x) = x^{-\lambda} T(y^{\lambda}f)(x), \qquad f \in C_c^{\infty}(\mathbb{R}_+),
$$
is also a local vector-valued Calder\'on-Zygmund operator, hence, given any $\delta \in \mathbb{R}$,
it extends to a bounded operator from $L^p(\mathbb{R}_+,x^{\delta}dx)$ into 
$L_B^p(\mathbb{R}_+,x^{\delta}dx)$, $1<p<\infty$, and from $L^1(\mathbb{R}_+,x^{\delta}dx)$
into $L^{1,\infty}_B(\mathbb{R}_+,x^{\delta}dx)$.
\end{lema}

This lemma can be justified, in a straightforward manner, by applying a vector-valued
variant of \cite[Theorem 4.3]{NS} and using the fact that each $A^p_{\textrm{loc}}$ class
(considered in \cite{NS}), $1\le p < \infty$, contains all power weights $x^{\delta}$, 
$\delta \in \mathbb{R}$.

Throughout the paper we use the convention that constants
may change their value (but not the dependence) from one occurrence to the next.
The notation $c_p$ means that the constant depends \emph{only} on $p$. Constants 
are always strictly positive and finite. Moreover, we distinguish ``big'' and ``small'' 
constants by using capital and small letters, respectively.

Finally, we shall implicitly use the simple fact that
$\sup_{t>0}t^\beta \exp(-\gamma t) = C_{\beta,\gamma} < \infty$ for arbitrary $\beta,\gamma >0$.

\section{The heat integral maximal operator} \label{sec:heat}

In this section we prove Theorem \ref{tmax}. Recall that for  $\lambda > -1\slash 2$
$$
W^{\lambda}_*f(x) = \sup_{t>0} \Big| \int_0^{\infty} W_t^{\lambda}(x,y)f(y)\, 
d\mu_{\lambda}(y)\Big|, \qquad x >0,
$$
with the heat kernel given by \eqref{heat_ker}.
We shall use the following estimates of $W_t^{\lambda}(x,y)$.

\begin{lema} \label{le21} 
Let $\lambda>-1/2$. Then for all $t,x,y>0$
$$
W_t^{\lambda}(x,y) \le C_{\lambda} \left\{
\begin{array}{ll}
x^{-2\lambda-1}, & y \le x\slash 2 \\
x^{-2\lambda-1} + t^{-1\slash 2}(xy)^{-\lambda} e^{-(x-y)^2\slash 4t}, & x\slash 2<y<2x \\
y^{-2\lambda-1} ({y^2}\slash {t})^{\lambda+1\slash 2}e^{-cy^2/t}, & 2x \le y
\end{array} \right. .
$$
\end{lema}

\begin{proof}
First observe that if $xy \le t$ then (\ref{P3}) implies
\begin{equation} \label{rr1}
W_t^\lambda(x,y) \le C_{\lambda} \frac{1}{{t}^{\lambda+1/2}} e^{-(x^2+y^2)\slash {4t}}
\le C_{\lambda} \frac{1}{y^{2\lambda+1}} \Big(\frac{y^2}{t}\Big)^{\lambda+1\slash 2} 
e^{-{y^2}\slash {4t}}.
\end{equation}
On the other hand, if $xy > t$ then (\ref{P4}) leads to
\begin{equation} \label{rr2}
W^{\lambda}_t(x,y) \le C_{\lambda}
\frac{1}{(xy)^{\lambda}\sqrt{t}} \, e^{-{(x-y)^2}\slash{4t}}.
\end{equation}

Assume now that $2x \le y$. Then, in view of \eqref{rr1}, it is enough to consider $xy>t$
and by \eqref{rr2} we get
\begin{align*}
W^{\lambda}_t(x,y)
\le  \frac{C_{\lambda}}{(xy)^{\lambda+1/2}} \Big(\frac{xy}{t}\Big)^{1\slash 2}
    e^{-{y^2}\slash{16t}}
 \le \frac{C_{\lambda}}{t^{\lambda+1/2}} \Big(\frac{y^2}{t}\Big)^{1\slash 2}
    e^{-{y^2}\slash{16t}}   =
\frac{C_{\lambda}}{y^{2\lambda+1}} \Big(\frac{y^2}{t}\Big)^{\lambda + 1}
    e^{-{y^2}\slash{16t}}.
\end{align*}
This implies the desired bound for $W_t^{\lambda}(x,y)$.

The case $y\le x\slash 2$, by the symmetry $W^\lambda_t(x,y) = W^\lambda_t(y,x)$,
is easily covered by a direct weakening of the already justified estimate for $2x\le y$.
Finally, the remaining estimate for comparable $x$ and $y$
follows by combining \eqref{rr1} and \eqref{rr2}.
\end{proof}

Consider the auxiliary maximal operator
$$
T f(x)=\sup_{t>0}\Big|\int_x^\infty \Big(\frac{y^2}{t}\Big)^{\lambda+1\slash 2} 
\exp\Big(-c\frac{y^2}{t}\Big) \frac{f(y)}{y} \,dy\Big|, \qquad x>0.
$$
Observe that $T$ can be controlled, up to a multiplicative constant, by the Hardy operator
$H^0_{\infty}$. Thus $T$ has all the mapping properties stated in Lemma \ref{Hardyinfty}
with $\eta =0$. Moreover, a straightforward computation shows that $T$ is bounded on
$L^{\infty}(\mathbb{R}_+)$. Further results concerning this and the more general operator
$$
T_\psi^\eta f(x)=\sup_{s}\Big|x^\eta\int_x^\infty f(y)\psi(s,y)y^{-\eta-1}\,dy\Big|, \qquad x>0,
$$
can be found in \cite[Lemma 3.3]{CH}, but they will not be used here except for Remark \ref{Re1} below.

We also invoke the local maximal function $M_{\textrm{loc}}^k$ defined by
$$
M_{\textrm{loc}}^k f(x)=\sup_{0<u<x<v<ku}\frac{1}{v-u}\int_u^v|f(y)|\, dy,\qquad x>0,
$$
for a given $k>1$. This operator, for any $\delta \in \mathbb{R}$, is bounded on
$L^p(\mathbb{R}_+,x^{\delta}dx)$, $1<p\le \infty$, and from $L^1(\mathbb{R}_+,x^{\delta}dx)$ 
to $L^{1,\infty}(\mathbb{R}_+,x^{\delta}dx)$, see \cite[Section 6]{NS}.

\subsection*{Proof of Theorem \ref{tmax}}
In order to show sufficiency parts in Theorem \ref{tmax} we split the kernel
$W^{\lambda}_t(x,y)$ according to the regions $0<y\le x\slash 2$, $x\slash 2 < y < 2x$,
$2x < y$, and denote the resulting maximal operators by
$N^{\lambda}_1$, $N^{\lambda}_2$ and $N^{\lambda}_3$, respectively. Then
$$
W^{\lambda}_*f(x) \le N^\lambda_1 |f|(x) + N^\lambda_2 |f|(x) + N^\lambda_3 |f|(x), \qquad x>0.
$$
Using Lemma \ref{le21} we get
$$
N^\lambda_1 |f|(x) \le C_{\lambda} H_0^{2\lambda} |f|(x), \qquad x>0.
$$
Another application of Lemma \ref{le21} gives
$$
N^\lambda_3 |f|(x) \le C_{\lambda} T |f|(x), \qquad x>0.
$$
Considering $N_2^{\lambda}$, again by Lemma \ref{le21} we have, for $x>0$,
\begin{align*}
N^{\lambda}_2 |f|(x) & \le C_{\lambda} \bigg(\int_{x/2}^{2x}|f(y)|\,
 \frac{dy}{y} + \sup_{t>0}\int_{x/2}^{2x}\frac{1}{\sqrt{t}}e^{-{(x-y)^2}\slash{4t}}
 |f(y)|\, dy\bigg) \le C_{\lambda} M^4_{\textrm{loc}} f(x).
\end{align*}

Now, taking into account the above estimates and facts, and appealing to Lemmas
\ref{Hardy0} and  \ref{Hardyinfty},  
we conclude the following mapping properties of the operator $W_*^{\lambda}$,
considered on the space $(\mathbb{R}_+,x^{\delta}dx)$. For
$1<p<\infty$, $W_*^{\lambda}$ is of strong type $(p,p)$ provided
that $-1<\delta< (2\lambda+1)p-1$. Moreover, $W_*^{\lambda}$ is of
strong type $(\infty,\infty)$ for each $\delta\in \mathbb{R}$. If
$-1<\delta\le 2\lambda$ then $W_*^{\lambda}$ is of weak type
$(1,1)$. Eventually, if $1<p<\infty$ and $-1 < \delta \le
(2\lambda+1)p-1$ then $W_*^{\lambda}$ is of restricted weak type
$(p,p)$. These facts, all together, justify the sufficiency parts in Theorem \ref{tmax}.

We pass to the proof of the necessity parts, that is  showing the
sharpness of the above results. Our task will be done once we
establish the following three statements (as before, we assume that
$\lambda>-1\slash 2$, $1\le p < \infty$ and the underlying space is $(\mathbb{R}_+,x^{\delta}dx)$).
\begin{itemize}
\item[(A)]
If $W_*^{\lambda}$ is of restricted weak type $(p,p)$ then $-1<\delta\le (2\lambda+1)p-1$.
\item[(B)]
$W_*^{\lambda}$ is not of weak type $(p,p)$ when $p>1$ and $\delta = (2\lambda+1)p-1$.
\item[(C)]
$W_*^{\lambda}$ is not of strong type $(1,1)$ if $-1 < \delta \le 2\lambda$.
\end{itemize}

To this end, let $f$ be a nonnegative function on $(0,\infty)$. Since by (\ref{P3})
$$
W_t^\lambda(x,y)\ge \frac{c_{\lambda}}{t^{\lambda+1/2}} e^{-{(x^2+y^2)}\slash{4t}}, \qquad xy<t,
$$
we see that
$$
W_t^{\lambda}f(x)\ge \frac{c_{\lambda}}{t^{\lambda+1/2}}\int_0^{t/x} e^{-{(x^2+y^2)}\slash {4t}}
    f(y) \, d\mu_{\lambda}(y), \qquad t,x>0.
$$
Consequently, since $W_*^{\lambda}f(x) \ge W_{x^2}^{\lambda}f(x)$, $x>0$, we get
\begin{equation} \label{W1}
W_*^{\lambda}f(x) 
\ge \frac{c_{\lambda}}{x^{2\lambda+1}}\int_0^x e^{-(x^2+y^2)\slash{4x^2}}f(y)\, d\mu_{\lambda}(y) 
\ge \frac{c_{\lambda}}{x^{2\lambda+1}}\int_0^x f(y)\,d\mu_{\lambda}(y). 
\end{equation}
Similarly, by (\ref{P4}) we have
\begin{equation} \label{WW3}
W_t^{\lambda}f(x)\ge \frac{c_{\lambda}}{\sqrt{t}}\int_{t/x}^\infty \frac{1}{(xy)^{\lambda}}
e^{-{(x-y)^2}\slash{4t}} f(y)\, d\mu_{\lambda}(y), \qquad t,x >0,
\end{equation}
and therefore
\begin{equation} \label{WW1}
W_*^{\lambda}f(x) \ge \frac{c_{\lambda}}{x}\int_{x}^\infty \frac{1}{(xy)^{\lambda}}
e^{-{(x-y)^2}\slash{4x^2}} f(y) \, d\mu_{\lambda}(y), \qquad x>0.
\end{equation}

Suppose now that $1\le p<\infty$ and $W_*^{\lambda}$ can be extended from
$L^2(\mathbb{R}_+,d\mu_{\lambda})$ to a restricted weak type $(p,p)$
operator on $(\mathbb{R}_+,x^\delta dx)$. Considering $f=\chi_{(1,2)}$, by \eqref{W1} we obtain
$$
W_*^{\lambda}f(x) \ge c_{\lambda} x^{-2\lambda-1}, \qquad x \ge 2.
$$
Then, by the weak type $(p,p)$ inequality satisfied by $f$, 
$$
\int_2^{\,\gamma^{-1\slash (2\lambda+1)}}y^\delta \, dy \le C_{p,\lambda} \, \gamma^{-p} 
$$
for $\gamma>0$ sufficiently small. It follows that the function
$\gamma \mapsto \gamma^{p-(\delta+1)/(2\lambda+1)}$ must be bounded for $\gamma$
near $0$ and we conclude that $\delta\le (2\lambda+1)p-1$.
On the other hand, in view of \eqref{WW1}, for $f=\chi_{(1,2)}$ as above and $x \in (0,1)$,
$$
W_*^{\lambda}f(x) \ge \frac{c_\lambda}{x^{\lambda+1}}\int_{1}^2
y^{\lambda} e^{-{(x-y)^2}\slash{4x^2}} dy \ge \frac{c_{\lambda}}{x^{\lambda+1}},
$$
which together with the weak type $(p,p)$ inequality for $f$ implies $\int_0^1 x^{\delta}dx < \infty$. 
The conclusion $\delta > -1$ follows and this completes proving statement (A).

Next, fix $1<p<\infty$, let $\delta=(2\lambda+1)p-1$ and suppose, on the contrary,
that $W_*^{\lambda}$ can be extended to a weak type $(p,p)$ operator on
$(\mathbb{R}_+,x^\delta dx)$. It is straightforward that with
$p'$ being the conjugate of $p$, $1\slash p + 1\slash p' =1$, the integral
$\int_0^1 x^{(2\lambda-\delta)p'+\delta}dx$ is infinite.
Therefore, for each $n \ge 1$ there exists a nonnegative function
$f_n\in L^2((0,1),d\mu_{\lambda})\cap L^p((0,1),x^\delta dx)$
such that $\|f_n\|_{L^p((0,1),x^\delta dx)}\le 1$ and
$\int_0^1 f_n(x)d\mu_{\lambda}(x)\ge n$. Extending $f_n$ to $\mathbb{R}_+$ by letting 
$f_n(x)=0$ for $x\ge 1$ and making use of \eqref{W1} we may write
$$
W_*^{\lambda}f_n(x)\ge c_{\lambda} n x^{-2\lambda-1}, \qquad x\ge 1, \quad n=1,2,\ldots.
$$
Now, in view of the weak type $(p,p)$ boundedness of $W_*^{\lambda}$, we get
$$
\int_1^{(n/(n-1))^{1/(2\lambda+1)}}y^\delta \, dy \le C_{\lambda}
    \bigg(\frac{\|f_n\|_{L^p((0,1),x^\delta dx)}}{n-1}\bigg)^p
        \le C_{\lambda} \frac{1}{(n-1)^{p}}, \qquad n=2,3,\ldots.
$$
This in turn implies boundedness of the sequence $\{n^p-(n-1)^p : n
\ge 2\}$, a contradiction because $p>1$. Thus statement (B) is justified.

It remains to prove that $W_*^{\lambda}$ is not strong type $(1,1)$ on
$(\mathbb{R}_+,x^\delta dx)$ when $-1<\delta\le 2\lambda$.
Assuming that $x\ge 1$ and $t \le 1\slash 2$ and restricting the interval of
integration in \eqref{WW3} we get
\begin{align*}
W_t^{\lambda}f(x) & \ge
c_{\lambda} \int_{x\slash 2}^{2x} \frac{1}{\sqrt{t}(xy)^{\lambda}}
e^{-{(x-y)^2}\slash{4t}} f(y) \,d\mu_{\lambda}(y) \\
& \ge c_{\lambda} \frac{1}{\sqrt{t}} \int_{x\slash 2}^{2x}
    e^{-{(x-y)^2}\slash{4t}} f(y)\,dy, \qquad x\ge 1, \quad t\le 1\slash 2.
\end{align*}
For $0<\varepsilon<1/2$ consider $1< y < 1+\varepsilon$ and $1+2\varepsilon < x < 2$.
Then obviously $x\slash 2<y<2x$ and, moreover, if $t=(x-1)^2\slash 2$ (notice that this quantity
is less than $1\slash 2$) we have
\begin{equation} \label{hasz}
\frac{(x-y)^2}{4t}\le \frac{1}{2}\Big(\frac{x-1+y-1}{x-1}\Big)^2 \le
    \frac{1}{2}\Big( 1 + \frac{\varepsilon}{2\varepsilon}\Big)^2 = \frac{9}{8}.
\end{equation}
Therefore, choosing $f_{\varepsilon}= \chi_{(1,1+\varepsilon)}$, we arrive at
$$
W_*^{\lambda}f_{\varepsilon}(x) \ge c_{\lambda} \frac{\varepsilon}{x-1},
    \qquad 1+2\varepsilon < x < 2.
$$
Consequently,
$$
\int_0^\infty |W_*^{\lambda}f_{\varepsilon}(x)|x^\delta \,dx
\ge c_{\lambda} \, \varepsilon \int_{1+2\varepsilon}^2\frac{x^\delta}{x-1}\,dx
\ge c_{\lambda,\delta} \, \varepsilon \int_{1+2\varepsilon}^2\frac{dx}{x-1}
=c_{\lambda,\delta}\, \varepsilon \log \frac{1}{2\varepsilon}.
$$
On the other hand, clearly $\|f_{\varepsilon}\|_{L^1(\mathbb{R}_+,x^\delta dx)}
\le C_{\delta}\, \varepsilon$. Letting $\varepsilon \to 0^+$ we see
that $W_*^{\lambda}$ is not bounded on $L^1(\mathbb{R}_+,x^\delta
dx)$. Statement (C) follows. 

The proof of Theorem \ref{tmax} is now complete.
\qed

\begin{remark} \label{Re1}
The boundedness properties of the maximal operator
$\widetilde{W}_*^\lambda$ can be obtained by following the lines of the
proof of Theorem \ref{tmax}, taking into account that
$\widetilde{W}_t^\lambda (x,y)=(xy)^\lambda W_t^\lambda(x,y)$.
In fact from Lemma \ref{le21} it follows that
$\widetilde{W}_*^\lambda f$ can be controlled by $H_0^\lambda
f+M^4_{\textrm{loc}}f+T_\psi^\lambda f$, with 
$\psi(s,y)=(y^2\slash s)^{\lambda+1\slash 2} e^{-cy^2\slash s}$. 
Invoking the mapping properties of $T^{\lambda}_{\psi}$ from \cite[Lemma 3.3]{CH},
and making use of \eqref{W1} and \eqref{WW1}, leads to similar results for
$\widetilde{W}^\lambda_*$ to those stated in Theorem \ref{tmax},
just replacing the interval $(-1,(2\lambda+1)p-1)$ by $(-\lambda p-1,(\lambda+1)p-1)$.
However, in the contrast with the $\Delta_{\lambda}$ setting,
these parallel results are not uniform in $\lambda$. A singularity occurs at $\lambda=0$,
exactly as described by the statements of \cite[Theorem 2.2]{CH} with $\alpha$
substituted by $\lambda-1\slash 2$.
\end{remark}

\section{Riesz transforms} \label{sec:riesz}

In this Section we prove Theorem \ref{triesz}. Recall that the Riesz transform
associated with the Bessel operator $\Delta_\lambda$ is formally defined by
$$
R_\lambda=D\Delta_\lambda^{-1/2}, \qquad \lambda>-1/2.
$$
We begin with a rigorous definition of the operator 
$\Delta_\lambda^{-1\slash 2}$. Then we obtain a representation of
$R_\lambda$ in terms of a principal value integral.

The negative power of $\Delta_\lambda$ can be defined, at least for
smooth functions with compact support, by
$$
\Delta_\lambda^{-1/2} f(x)=
\frac{1}{\sqrt{\pi}}\int_0^\infty \Big(W_t^{\lambda}f(x)- \chi_{\{\lambda \le 0\}}
    W_t^{\lambda} f(0) \Big) \frac{dt}{\sqrt{t}},
$$
with $W_t^{\lambda}f(0)$ understood as $\lim_{x\to 0^+} W_t^{\lambda}f(x)$.
Note that a compensating term is necessary
to make the integral convergent when $-1\slash 2 < \lambda \le 0$, as can be seen below.

\begin{propo}
Let $\lambda>-1\slash 2$ and $f \in C^{\infty}_c(\mathbb{R}_+)$.
Then the limit $W_t^{\lambda}f(0)$ exists for each fixed $t>0$,
and the function $\Delta_{\lambda}^{-1\slash 2} f(x)$, $x>0$, is well defined. 
\end{propo}

\begin{proof}
Suppose that $\supp f \subset (a,b)$, $0<a<b<\infty$. To see that
the limit exists it is sufficient to notice that, in view of \eqref{P3}, we have
$$
W_t^{\lambda}(0,y) := \lim_{x\to 0^+}W_t^{\lambda}(x,y) = \frac{1}{2^{2\lambda}\Gamma(\lambda+1\slash 2)}
\frac{1}{t^{\lambda+1\slash 2}} e^{-y^2\slash 4t},
$$
the convergence being uniform in $y<b$.

Now fix $x>0$. By using (\ref{P1}) we get
\begin{align*}
& \int_1^\infty \Big|W_t^{\lambda}f(x)- \chi_{\{\lambda \le 0\}}
    W_t^{\lambda} f(0) \Big| \frac{dt}{\sqrt{t}} \\
&\le \int_1^\infty \int_a^b \bigg|\frac{(xy)^{-\lambda+1/2}}{2t}
I_{\lambda-1/2}\Big(\frac{xy}{2t}\Big)e^{-(x^2+y^2)/4t}-
\frac{\chi_{\{\lambda \le 0\}}e^{-y^2\slash 4t}}{t^{\lambda+1/2}2^{2\lambda}\Gamma(\lambda+1/2)}
\bigg| |f(y)|\,d\mu_{\lambda}(y)\frac{dt}{\sqrt{t}}\\
&\le C_{\lambda,b,x} \int_1^\infty \int_a^b \Big(\big|e^{-(x^2+y^2)/4t}-
\chi_{\{\lambda \le 0\}}e^{-y^2\slash 4t}\big| + \frac{(xy)^2}{t^2} \Big) |f(y)|\, d\mu_{\lambda}(y)
\frac{dt}{t^{\lambda+1}}\\
&\le C_{\lambda,b,x}  \int_a^b \int_1^\infty
\Big(\frac{\chi_{\{\lambda \le 0\}} {x^2}}{t^{\lambda+2}} +
\frac{\chi_{\{\lambda >0\}}}{t^{\lambda+1}}  +
\frac{(xy)^2}{t^{\lambda+3}} \Big) \,dt \, |f(y)|\, d\mu_{\lambda}(y)
\end{align*}
and the last double integral converges. Further, an application of (\ref{P4}) gives
\begin{align*}
& \int_0^1 \Big|W_t^{\lambda}f(x)- \chi_{\{\lambda \le 0\}}
    W_t^{\lambda} f(0) \Big| \frac{dt}{\sqrt{t}} \\
& \le C_{\lambda,a,x} \int_0^1 \int_a^b \Big(\frac{(xy)^{-\lambda}}{\sqrt{t}}
    e^{-(x-y)^2/4t} + \frac{\chi_{\{\lambda\le 0\}}}{t^{\lambda+1\slash 2}} 
    e^{-y^2\slash 4t}\Big) |f(y)| \, d\mu_{\lambda}(y)\frac{dt}{\sqrt{t}}\\
&\le C_{\lambda,a,b,x} \int_a^b \bigg\{ \frac{1}{|x-y|^{1\slash 2}}
    \int_0^1 \Big(\frac{(x-y)^{2}}{t}\Big)^{1\slash 4}
    e^{-(x-y)^2/4t}\frac{dt}{t^{3\slash 4}} \\ & \quad + \chi_{\{\lambda \le 0\}}
\frac{1}{y^{1\slash 2}}\int_0^1 \Big(\frac{y^{2}}{t}\Big)^{1\slash 4}e^{-y^2/4t}
    \frac{dt}{t^{\lambda+3\slash 4}} \bigg\} |f(y)|\,dy .
\end{align*}
The last expression is controlled, up to a multiplicative constant, by the convergent integral
$\int_a^b (|x-y|^{-1\slash 2}+\chi_{\{\lambda\le 0\}} y^{-1\slash 2})|f(y)| dy$.
Combining the above facts we conclude that
the integral defining $\Delta_\lambda^{-1/2} f(x)$ converges absolutely.
\end{proof}

\begin{propo} \label{pvalue}
Let $\lambda > -1\slash 2$. For $f\in C_c^\infty(\mathbb{R}_+)$ the function
$\Delta_{\lambda}^{-1\slash 2}f$ is differentiable and
\begin{equation*} 
R_{\lambda}f(x) := D\Delta_\lambda^{-1/2}f(x) =
\pv \int_0^\infty R_{\lambda}(x,y) f(y)\,d\mu_{\lambda}(y), \qquad x>0,
\end{equation*}
with the kernel given by
$$
R_\lambda(x,y)=\frac{1}{\sqrt{\pi}}\int_0^\infty \frac{\partial}{\partial x}W_t^\lambda(x,y)
    \frac{dt}{\sqrt{t}}, \qquad x,y >0, \quad x \neq y.
$$
\end{propo}

\begin{proof}
Our reasoning is based on a comparison with the classical setting, which allows to control the
essential singularity. Let $\mathcal{H}$ denote the Hilbert transform and
let $f \in C_c^{\infty}(\mathbb{R})$. It is well known that
\begin{align}
\mathcal{H}f(x) & = \frac{1}{\pi}\pv \int_{\mathbb{R}} \frac{f(y)}{y-x} dy \nonumber \\
& = \pv \int_{\mathbb{R}} \bigg( \frac{1}{\sqrt{\pi}} \int_0^{\infty}
    \frac{\partial}{\partial x} \mathcal{W}_t(x,y) \frac{dt}{\sqrt{t}}\bigg)f(y)\,dy \nonumber \\
& = \frac{1}{\sqrt{\pi}} \frac{\partial}{\partial{x}}
    \int_0^{\infty} \Big( \mathcal{W}_t f(x) - \mathcal{W}_t f(0) \Big)
    \frac{dt}{\sqrt{t}}, \qquad x \in \mathbb{R}. \label{hilb}
\end{align}
Thus we decompose
\begin{align*}
& D \Delta_{\lambda}^{-1\slash 2} f(x) \\
& = \frac{1}{\sqrt{\pi}} x^{-\lambda} \frac{\partial}{\partial x}
    \int_0^{\infty} \Big( \mathcal{W}_t (y^{\lambda}f)(x) -
    \mathcal{W}_t (y^{\lambda}f)(0) \Big) \frac{dt}{\sqrt{t}} \\
    &\quad - \frac{\lambda}{\sqrt{\pi}} x^{-\lambda-1} \int_0^{\infty} \Big( \mathcal{W}_t
     (y^{\lambda}f)(x) - \mathcal{W}_t (y^{\lambda}f)(0) \Big) \frac{dt}{\sqrt{t}} \\
    & \quad + \frac{1}{\sqrt{\pi}} \frac{\partial}{\partial x}
        \int_0^{\infty} \Big( W_t^{\lambda}f(x) - \chi_{\{\lambda \le 0\}}
        W_t^{\lambda} f(0) - x^{-\lambda} \mathcal{W}_t(y^{\lambda}f)(x) + x^{-\lambda}
        \mathcal{W}_t(y^{\lambda}f)(0) \Big) \frac{dt}{\sqrt{t}}\\
    & \equiv \mathcal{I}_1 - \mathcal{I}_2 + \mathcal{I}_3.
\end{align*}
In view of \eqref{hilb} it follows that
$$
\mathcal{I}_1 = \pv \int_{0}^{\infty} \bigg( \frac{1}{\sqrt{\pi}} \int_0^{\infty}
    \frac{\partial}{\partial x} \mathcal{W}_t(x,y) \frac{dt}{\sqrt{t}}\bigg)(xy)^{-\lambda}
     f(y)\,d\mu_{\lambda}(y);
$$
this term contains the crucial singularity and, as we shall see, no
singular integrals emerge from $\mathcal{I}_2$ and $\mathcal{I}_3$.

Next, we prove that it is legitimate to pass with $\frac{\partial}{\partial x}$ under the integral
sign in $\mathcal{I}_3$.
This task, however, is directly reduced to showing that, for any fixed $0<a<b< \infty$, the quantities
\begin{align*}
\mathcal{J}_1(t) & = \sup_{x \in [a,b]} \big| \mathcal{W}_t (y^{\lambda}f)(x)
        - \mathcal{W}_t (y^{\lambda}f)(0) \big|,\\
\mathcal{J}_2(t) & = \sup_{x \in [a,b]} \Big| \frac{\partial}{\partial x} W_t^{\lambda}f(x)
    - x^{-\lambda} \frac{\partial}{\partial x} \mathcal{W}_t(y^{\lambda}f)(x) \Big|,
\end{align*}
can be majorized by functions of $t$ belonging to
$L^1(\mathbb{R}_+,t^{-1\slash 2}dt)$. In addition, without any
loss of generality, it may be assumed that $\supp f \subset
[a,b]$, with $0<a<b<\infty$, and $\|y^{\lambda}f\|_{\infty}\le 1$, $\|f\|_{\infty}\le 1$.

Observe that
$$
\mathcal{J}_1(t) \le \frac{1}{t^{1\slash 2}} \sup_{x \in [a,b]} \int_a^b \big| e^{-{(x-y)^2}\slash{4t}}
    - e^{-{y^2}\slash{4t}} \big| dy.
$$
For $t \ge 1$ and $x,y \in [a,b]$ we have
$$
\big| e^{-{(x-y)^2}\slash{4t}} - e^{-{y^2}\slash{4t}} \big| = \bigg| \frac{y^2}{4t} -
    \frac{(x-y)^2}{4t} + \mathcal{O} \Big(\frac{(x-y)^4}{t^2} \Big) -
    \mathcal{O}\Big(\frac{y^4}{t^2}\Big) \bigg| \le \frac{C_b}{t}.
$$
Also, since for all $t>0$ and $x,y \in [a,b]$,
$$
\big| e^{-{(x-y)^2}\slash{4t}} - e^{-{y^2}\slash{4t}} \big| \le
    e^{-{(x-y)^2}\slash{4t}} + e^{-{a^2}\slash{4t}} ,
$$
we get
$$
\frac{1}{t^{1\slash 2}} \int_a^b \big| e^{-{(x-y)^2}\slash{4t}}
    - e^{-{y^2}\slash{4t}} \big| dy\le \frac{2}{t^{1\slash 2}}\int_{-\infty}^\infty e^{-u^2/4t}du
    = 4\sqrt{\pi}.
$$
Then
$$
\mathcal{J}_1(t) \le C_{b} \min\big( t^{-3\slash 2}, 1\big) \in
    L^1(\mathbb{R_+},t^{-1\slash 2}dt).
$$

The treatment of $\mathcal{J}_2$ is not as straightforward. We consider two cases.

\noindent \textit{Case 1}: $t \le b^2$.
By \eqref{P2} we get
\begin{align*}
\frac{\partial}{\partial x}W_t^\lambda(x,y) & = \frac{1}{(2t)^\lambda}e^{-xy/2t}
\frac{\partial}{\partial x}\Big(\frac{1}{\sqrt{2t}}e^{-(x-y)^2/4t}\Big)
\Big(\frac{xy}{2t}\Big)^{-\lambda+1/2} I_{\lambda-1/2}\Big(\frac{xy}{2t}\Big)\\
&\quad -\frac{y}{(2t)^{\lambda+3/2}}\Big(\frac{xy}{2t}\Big)^{-\lambda+1\slash 2}
\bigg( I_{\lambda-1/2}\Big(\frac{xy}{2t}\Big) 
-  I_{\lambda+1/2}\Big(\frac{xy}{2t}\Big)\bigg) e^{-(x^2+y^2)/4t}\\
&\equiv H_{\lambda,1}(t,x,y) - H_{\lambda,2}(t,x,y), \qquad t,x,y>0.
\end{align*}
Next observe that \eqref{P4} implies
\begin{align}
H_{\lambda,1}(t,x,y)&= \frac{1}{\sqrt{2\pi}}\frac{\partial}{\partial x}
\Big(\frac{1}{\sqrt{2t}}e^{-(x-y)^2/4t}\Big)(xy)^{-\lambda}
\bigg(1+\mathcal{O}\Big(\frac{t}{xy}\Big)\bigg) \nonumber \\
&= (xy)^{-\lambda} \frac{\partial}{\partial x} \mathcal{W}_t(x,y)
-\frac{(xy)^{-\lambda}}{t^{3\slash 2}} (x-y) e^{-(x-y)^2/4t}
\mathcal{O}\Big(\frac{t}{xy}\Big), \qquad x,y\in [a,b]. \label{T_1}
\end{align}
In order to analyze $H_{\lambda,2}(t,x,y)$ we use again \eqref{P4} and obtain
\begin{equation} \label{T_2}
|H_{\lambda,2}(t,x,y)|=\frac{y}{t^{\lambda+3/2}}e^{-(x-y)^2/4t}\Big(\frac{xy}{2t}\Big)^{-\lambda}
\mathcal{O}\Big(\frac{t}{xy}\Big),\qquad x,y\in [a,b].
\end{equation}
Thus $H_{\lambda,1}$ and $H_{\lambda,2}$ satisfy
\begin{align*}
\Big| H_{\lambda,1}(t,x,y) - (xy)^{-\lambda} \frac{\partial}{\partial x} \mathcal{W}_t(x,y) \Big| & \le
    C_{\lambda,a,b}  \frac{1}{{t^{1\slash 2}}} e^{-(x-y)^2\slash 4t},\qquad x,y\in [a,b],\\
|H_{\lambda,2}(t,x,y)| & \le C_{\lambda,a,b}  \frac{1}{t^{1\slash
2}} e^{-(x-y)^2\slash 4t}, \qquad x,y\in [a,b].
\end{align*}
Consequently,
$$
\int_a^b\Big| H_{\lambda,1}(t,x,y) - (xy)^{-\lambda}
\frac{\partial}{\partial x} \mathcal{W}_t(x,y)
\Big||f(y)|\,d\mu_{\lambda}(y)\le C_{\lambda,a,b}\frac{1}{t^{1\slash
2}}\int_{-\infty}^\infty e^{-(x-y)^2\slash 4t}dy = C_{\lambda,a,b},
$$
and
$$
\int_a^b|H_{\lambda,2}(t,x,y)||f(y)|\,d\mu_{\lambda}(y) = C_{\lambda,a,b}.
$$

Hence $\mathcal{J}_2(t) \le C_{\lambda,a,b}$ when $t \le b^2$.

\noindent \textit{Case 2}: $t>a^2$.
We deduce from \eqref{P2} that
\begin{align}
\frac{\partial}{\partial x}W_t^\lambda(x,y)
& =\frac{1}{(2t)^{\lambda+1/2}}e^{-(x^2+y^2)/4t}\bigg(x\Big(\frac{y}{2t}\Big)^2
\Big(\frac{xy}{2t}\Big)^{-\lambda-1/2}I_{\lambda+1/2}\Big(\frac{xy}{2t}\Big)\nonumber\\\label{45}
&\quad -\frac{x}{2t}\Big(\frac{xy}{2t}\Big)^{-\lambda+1/2}I_{\lambda-1/2}
\Big(\frac{xy}{2t}\Big)\bigg).
\end{align}
Then \eqref{P3} implies
\begin{equation} \label{T_3}
\Big|\frac{\partial}{\partial x}W_t^\lambda(x,y)\Big|\le C_{\lambda,a,b}
\frac{x}{t^{\lambda+3/2}}e^{-(x^2+y^2)/4t}\Big(\frac{y^2}{t}+1\Big)\le
\frac{C_{\lambda,a,b}}{t^{\lambda+3/2}},\qquad x,y\in [a,b].
\end{equation}
Therefore
$$
\int_a^b \Big| \frac{\partial}{\partial x} W_t^{\lambda}(x,y) f(y) \Big| d\mu_{\lambda}(y) \le
    \frac{C_{\lambda,a,b}}{t^{\lambda+3\slash 2}},\qquad x\in [a,b].
$$
On the other hand, for every $x\in [a,b]$,
$$
\int_a^b (xy)^{-\lambda} \Big| \frac{\partial}{\partial x} \mathcal{W}_t(x,y) f(y) \Big|
    d\mu_{\lambda}(y) \le \int_a^b (xy)^{-\lambda} \frac{|x-y|}{t^{3\slash 2}} 
    e^{-(x-y)^2\slash 4t} |f(y)| d\mu_{\lambda}(y) \le \frac{C_{\lambda,a,b}}{t^{3\slash 2}}.
$$
In this way we conclude that $\mathcal{J}_2(t) \le C_{\lambda,a,b} t^{-1}$ for $t>a^2$.

A combination of Case 1 and Case 2 reveals that
$$
\mathcal{J}_2(t) \le C_{\lambda,a,b} \min\big( 1,t^{-1} \big)
    \in L^1(\mathbb{R}_+,t^{-1\slash 2}dt).
$$

Now passing with $\frac{\partial}{\partial x}$ under the
integral in $\mathcal{I}_3$ is justified and we get
\begin{align*}
D \Delta_{\lambda}^{-1\slash 2} f(x) & = \pv \int_{0}^{\infty} \bigg( \frac{1}{\sqrt{\pi}}
 \int_0^{\infty} \frac{\partial}{\partial x} \mathcal{W}_t(x,y) \frac{dt}{\sqrt{t}}\bigg)(xy)^{-\lambda}
  f(y)\,d\mu_{\lambda}(y)\\
 & \quad + \frac{1}{\sqrt{\pi}} \int_0^{\infty} \Big( \frac{\partial}{\partial x} W_t^{\lambda}f(x) -
  x^{-\lambda} \frac{\partial}{\partial x} \mathcal{W}_t (y^{\lambda}f)(x) \Big) \frac{dt}{\sqrt{t}} \\
 & \equiv \mathcal{I}_4 + \mathcal{I}_5.
\end{align*}
Having in mind the estimates obtained so far in this proof, it is straightforward to check that
$$
\mathcal{I}_5 = \frac{1}{\sqrt{\pi}} \int_0^{\infty} \int_0^{\infty} \Big(
    \frac{\partial}{\partial x} W_t^{\lambda}(x,y) - (xy)^{-\lambda} \frac{\partial}{\partial x}
    \mathcal{W}_t(x,y) \Big) f(y) \, d\mu_{\lambda}(y) \frac{dt}{\sqrt{t}}
$$
and that the double integral above converges absolutely for any fixed $x>0$. Thus the order of
integration may be switched and then the cancellations occurring between
$\mathcal{I}_4$ and $\mathcal{I}_5$ lead to the desired principal value integral representation of
$D \Delta_{\lambda}^{-1\slash 2}f(x)$. Notice that the differentiability of 
$\Delta_{\lambda}^{-1\slash 2}f$ is implicitly contained in the whole reasoning.
\end{proof}

It is perhaps worth to mention that the kernel $R_{\lambda}(x,y)$ can be expressed explicitly
in terms of the Gauss hypergeometric function $_2F_1$
(for the definition see for instance \cite[Chapter 9]{Leb}). More precisely, by means of
\eqref{P2} and the integral formula, cf. \cite[2.15.3 (2)]{prud},
\begin{equation} \label{intF}
\int_0^{\infty} z^{1\slash 2} e^{-\alpha z}I_{\nu}(\beta z)\, dz = \alpha^{-\nu-3\slash 2}
    \Big(\frac{\beta}{2}\Big)^{\nu} \frac{\Gamma(\nu+3\slash 2)}{\Gamma(\nu+1)}\;
        {_2F_1} \bigg( \frac{\nu+3\slash 2}{2}, \frac{\nu+5\slash 2}{2}; \nu+1;
            \frac{\beta^2}{\alpha^2} \bigg),
\end{equation}
valid for $\nu>-3\slash 2$ and $\alpha>\beta>0$, one computes that for $x \neq y$
\begin{align*}
R_{\lambda}(x,y) = & \frac{2}{\sqrt{\pi}} \frac{\Gamma(\lambda + 2)}{\Gamma(\lambda+3\slash 2)}
    (xy)^{-\lambda-1}
    \bigg\{ y \Phi_{x,y}^{\lambda+2}
     \;{_2F_1}\Big( \frac{\lambda + 2}{2},
        \frac{\lambda+3}{2}; \frac{2\lambda+3}{2}; 4 \Phi_{x,y}^2\Big)\\
& - \frac{\lambda+1\slash 2}{\lambda+1} \;x \Phi_{x,y}^{\lambda+1}
     \;{_2F_1}\Big( \frac{\lambda+1}{2},
        \frac{\lambda+2}{2}; \frac{2\lambda+1}{2}; 4 \Phi_{x,y}^2\Big)  \bigg\},
\end{align*}
with $\Phi_{x,y} = xy\slash (x^2+y^2)$. This representation, even though explicit,
does not seem to be convenient for performing necessary kernel estimates since there
are essential cancellations between the two terms containing $_2F_1$ functions with
different parameters. For $\lambda=0$ the above expression can be simplified (the relevant
property of $_2F_1$ can be found in \cite[Section 9.8]{Leb}) and we have
\begin{equation} \label{rkzero}
R_0(x,y) = \frac{1}{\pi} \Big( \frac{1}{y-x} - \frac{1}{y+x}
\Big), \qquad x \neq y.
\end{equation}
Note that the same result can be obtained more directly since
$W_t^0(x,y) = \mathcal{W}_t(x,y)+\mathcal{W}_t(x,-y)$ (this identity follows from the fact that
$I_{-1\slash 2}(z) = \sqrt{2\slash \pi z} \cosh z$, cf. \cite[(5.8.5)]{Leb}) and the conclusion
is obtained by a comparison with the Hilbert transform kernel.
Finally, by \eqref{rkzero} we see that $R_0 f$ coincides with the Hilbert transform
of the even extension of $f$, restricted to the positive half-line.

The following estimates for the kernel $R_\lambda (x,y)$ will be
crucial in proving Theorem \ref{triesz}. They can be
obtained also as consequences of \cite[(1.6)]{Anker}, but
our procedure, via the heat kernel, is different from that 
contained in \cite[Lemma 2.1]{Ker1}, \cite[Theorem 2.1]{Ker2} and proving the estimate in \cite{Anker}. 

\begin{lema} \label{leriesz}
Let $\lambda>-1/2$. Then for all $x,y>0$, $x\neq y$, the integral defining $R_{\lambda}(x,y)$
converges absolutely and we have
$$
R_\lambda(x,y)=\frac{1}{{\pi}}\frac{(xy)^{-\lambda}}{y-x}+\mathcal{O}
    \bigg(y^{-2\lambda-1}\Big(1+\log\frac{xy}{(y-x)^2}\Big)\bigg), \qquad x/2<y<2x.
$$
Moreover, in the off-diagonal region,
$$
|R_{\lambda}(x,y)| \le C_{\lambda} \left\{
\begin{array}{ll}
x^{-2\lambda-1}, & y \le x\slash 2 \\
x y^{-2\lambda-2}, & 2x \le y
\end{array} \right. .
$$
\end{lema}

\begin{proof}
Observe first that by (\ref{P2})
\begin{align}
\frac{\partial}{\partial x} W_t^{\lambda}(x,y) = &
\frac{1}{(2t)^{\lambda+1/2}}\bigg(x \Big(\frac{y}{2t}\Big)^2
\Big(\frac{xy}{2t}\Big)^{-\lambda-1/2} I_{\lambda+1/2}\Big(\frac{xy}{2t}\Big)\nonumber\\
&-\frac{x}{2t}\Big(\frac{xy}{2t}\Big)^{-\lambda+1/2}
I_{\lambda-1/2}\Big(\frac{xy}{2t}\Big)\bigg) e^{-(x^2+y^2)/4t}, \qquad t,x,y>0. \label{C1}
\end{align}
Now assume that $0<x/2<y<2x$ and $x \neq y$. We write
$$
\sqrt{\pi} R_\lambda(x,y) = \int_0^{xy} \frac{\partial}{\partial x}W_t^\lambda(x,y)\frac{dt}{\sqrt{t}}
+ \int_{xy}^\infty \frac{\partial}{\partial x}W_t^\lambda(x,y)\frac{dt}{\sqrt{t}}
\equiv R_\lambda^1(x,y)+R_\lambda^2(x,y).
$$
By (\ref{C1}) and (\ref{P3}) it follows that
$$
|R_\lambda^2(x,y)| \le C_{\lambda} \int_{xy}^\infty
\Big(\frac{xy^2}{t^2}+\frac{x}{t}\Big)\frac{dt}{t^{\lambda+1}}
\le C_{\lambda} \Big(\frac{xy^2}{(xy)^{\lambda+2}}+\frac{x}{(xy)^{\lambda+1}}\Big)
\le C_{\lambda} y^{-2\lambda-1}.
$$
Next, by (\ref{C1}) and (\ref{P4}) we get
\begin{align*}
&R_\lambda^1(x,y)\\
&= \int_0^{xy}
\frac{2^{-\lambda-1\slash 2}}{ t^{\lambda+1}}\Big(\frac{xy}{2t}\Big)^{-\lambda} \bigg( \frac{y}{2t}\Big(\frac{xy}{2t}\Big)^{1/2} I_{\lambda+1/2}\Big(\frac{xy}{2t}\Big) - \frac{x}{2t}\Big(\frac{xy}{2t}\Big)^{1/2}I_{\lambda-1/2}\Big(\frac{xy}{2t}\Big)\bigg)
e^{-(x^2+y^2)/4t}dt\\
& =\frac{1}{2\sqrt{2}}\int_0^{xy} \frac{(xy)^{-\lambda}}{t^2}\bigg(\frac{y-x}{\sqrt{2\pi}}+
\mathcal{O}\Big(\frac{t}{y}\Big)\bigg) e^{-(x-y)^2/4t}dt\\
& =\frac{1}{4\sqrt{\pi}}(xy)^{-\lambda}(y-x)\int_0^\infty \frac{1}{t^2} e^{-(x-y)^2/4t}dt- \frac{1}{4\sqrt{\pi}}(xy)^{-\lambda}(y-x)\int_{xy}^\infty \frac{1}{t^2}
e^{-(x-y)^2/4t}dt\\
& \quad + (xy)^{-\lambda}\int_0^{xy}\frac{1}{t^2}\mathcal{O}\Big(\frac{t}{y}\Big)e^{-(x-y)^2/4t}dt\\
& \equiv K^{\lambda}_1(x,y)+K^{\lambda}_2(x,y)+K^{\lambda}_3(x,y).
\end{align*}
We analyze each $K^{\lambda}_i(x,y)$, $i=1,2,3$, separately. By a direct computation
$$
K^{\lambda}_1(x,y) = 
\frac{1}{\sqrt{\pi}}\frac{(xy)^{-\lambda}}{y-x}.
$$
To estimate $K^{\lambda}_2(x,y)$ we write
$$
|K^{\lambda}_2(x,y)| \le (xy)^{-\lambda}|y-x|\int_{xy}^\infty
\frac{1}{t^2}e^{-(x-y)^2/4t}dt
\le (xy)^{-\lambda}x\int_{xy}^\infty \frac{dt}{t^2}
\le C_{\lambda} y^{-2\lambda-1}.
$$
In case of $K^{\lambda}_3(x,y)$ one has
$$
|K^{\lambda}_3(x,y)| \le C_{\lambda} (xy)^{-\lambda}\frac{1}{y}\int_0^{xy} \frac{1}{t}e^{-(x-y)^2/4t}dt
\le C_{\lambda} y^{-2\lambda-1}\int_{(x-y)^2/4xy}^\infty
\frac{1}{u} e^{-u}du,
$$
and splitting the last integral according to $u<1$ and $u\ge 1$ it becomes clear that
$$
|K^{\lambda}_3(x,y)| \le C_{\lambda} y^{-2\lambda-1}\Big(1+\log\frac{xy}{(x-y)^2}\Big).
$$
This completes the proof of the diagonal estimate of the lemma.

In order to justify the remaining estimates it suffices to bound suitably 
$R_{\lambda}^i(x,y)$, $i=1,2$, in the off-diagonal region.
Notice that (\ref{P3}), together with \eqref{C1}, implies
\begin{align*}
& |R_\lambda^2(x,y)| \\
& \le C_{\lambda}x\int_{xy}^\infty
\frac{1}{t^{\lambda+2}}\Big(\frac{y^2}{t}+1\Big)e^{-(x^2+y^2)/4t}dt\\
&\le C_{\lambda}x\bigg(\frac{y^2}{(x^2+y^2)^{\lambda+2}}\int_0^{(x^2+y^2)/4xy} u^{\lambda+1}e^{-u}du
+ \frac{1}{(x^2+y^2)^{\lambda+1}}\int_0^{(x^2+y^2)/4xy}
u^{\lambda}e^{-u}du\bigg)\\
&\le C_{\lambda} \frac{x}{(x^2+y^2)^{\lambda+1}}.
\end{align*}
From this we easily obtain the desired estimates, with
$R_{\lambda}(x,y)$ replaced by $R^2_{\lambda}(x,y)$.
Finally, by applying (\ref{P4}) together with \eqref{C1},
$$
|R_\lambda^1(x,y)| \le C_{\lambda} (xy)^{-\lambda}(x+y)\int_0^{xy}
e^{-(x-y)^2/4t}\frac{dt}{t^2} \le C_{\lambda}
\frac{(xy)^{-\lambda}}{|x-y|^{2\lambda+3}}\int_0^{xy}t^\lambda dt
\le C_{\lambda} \frac{xy}{|x-y|^{2\lambda+3}}
$$
and the off-diagonal estimates follow again, this time for $R_{\lambda}^1(x,y)$.
The proof of Lemma~\ref{leriesz} is finished.
\end{proof}

We now show that the off-diagonal estimates of Lemma \ref{leriesz} are sharp in certain regions.

\begin{lema} \label{l2riesz}
Let $\lambda>-1/2$. There exist $b>1$ and a (positive) constant $c_{\lambda}$ such that
\begin{align*}
R_{\lambda}(x,y) & \le - c_{\lambda} x^{-2\lambda-1}, \qquad 0 < y \le x\slash b,\\
R_{\lambda}(x,y) & \ge c_{\lambda} x y^{-2\lambda-2}, \qquad 0 < bx \le y.
\end{align*}
\end{lema}

\begin{proof}
Using \eqref{C1} and then performing the change of variable $u=2t\slash x^2$ we obtain
\begin{align*}
R_\lambda(x,y) & =  \frac{1}{2^{\lambda+1/2}}\int_0^\infty
\bigg(x \Big(\frac{y}{2t}\Big)^2
\Big(\frac{xy}{2t}\Big)^{-\lambda-1/2} I_{\lambda+1/2}\Big(\frac{xy}{2t}\Big) \\
& \quad -\frac{x}{2t}\Big(\frac{xy}{2t}\Big)^{-\lambda+1/2}
I_{\lambda-1/2}\Big(\frac{xy}{2t}\Big)\bigg)e^{-(x^2+y^2)/4t}\frac{dt}{t^{\lambda+1}}\\
& = \frac{1}{\sqrt{2}x^{2\lambda+1}}\int_0^\infty
\frac{1}{u^{\lambda+1}}\bigg(\Big(\frac{z}{u}\Big)^2\Big(\frac{z}{u}\Big)^{-\lambda-1/2}
I_{\lambda+1/2}\Big(\frac{z}{u}\Big)\\
& \quad -\frac{1}{u}\Big(\frac{z}{u}\Big)^{-\lambda+1/2}I_{\lambda-1/2}\Big(\frac{z}{u}\Big)\bigg)
e^{-(1+z^2)/2u} du,\qquad x,y>0,
\end{align*}
where $z=y/x$. Then, with the aid of \eqref{P3} and the dominated convergence theorem
(a suitable integrable majorant can be derived by means of \eqref{P3} and \eqref{P4}), we see that
$$
\lim_{z\to 0^{+}}x^{2\lambda+1}R_\lambda(x,y) = - \frac{1}{2^{\lambda}\Gamma(\lambda+1/2)}
\int_0^\infty \frac{1}{u^{\lambda+2}}e^{-1/2u} du
= -  \frac{2\,\Gamma(\lambda+1)}{\Gamma(\lambda+1/2)},
$$
and the desired bound for $y < x\slash b$ follows.
Similarly, changing the variable $u = 2t\slash y^2$, we can write
\begin{align*}
R_\lambda(x,y) & = \frac{x}{\sqrt{2} y^{2\lambda+2}}\int_0^\infty
\frac{1}{u^{\lambda+1}}\bigg(\frac{1}{u^2}\Big(\frac{z}{u}\Big)^{-\lambda-1/2}
I_{\lambda+1/2}\Big(\frac{z}{u}\Big)\\
&\quad -\frac{1}{u} \Big(\frac{z}{u}\Big)^{-\lambda+1/2}I_{\lambda-1/2}
\Big(\frac{z}{u}\Big)\bigg) e^{-(1+z^2)/2u}du, \qquad x,y>0,
\end{align*}
being now $z=x/y$, thus, again in view of \eqref{P3} and the dominated convergence theorem,
$$ 
\lim_{z\to 0^+}\frac{y^{2\lambda+2}}{x}R_\lambda(x,y)
=\frac{1}{2^{\lambda+1}\Gamma(\lambda+1/2)}
\int_0^\infty\frac{1}{u^{\lambda+2}}\Big(\frac{1}{(\lambda+1/2)u}-2\Big)e^{-1/2u}du
=\frac{\Gamma(\lambda+1)}{\Gamma(\lambda+3/2)}.
$$ 
This gives the remaining bound for $bx<y$.
\end{proof}

Consider the auxiliary operator
$$
\mathcal{H}_{\lambda, \textrm{loc}}f(x) = \frac{1}{\pi} \pv
\int_{x/2}^{2x}\frac{(xy)^{-\lambda}}{y-x}f(y)\,d\mu_{\lambda}(y), \qquad x>0,
$$
defined for, say, $f \in C_c^{\infty}(\mathbb{R}_+)$.
Note that the principal value integral converges for a.e. $x>0$ if only
$f$ is locally integrable; this follows by the relation with the Hilbert transform.

According to \cite[Lemma 1]{AM}, the local Hilbert transform $\mathcal{H}_{0, \textrm{loc}}$ 
is bounded on $L^2(\mathbb{R}_+,dx)$ (this fact can be also proved directly by combining the
$L^2$-boundedness of the Hilbert transform with classic Hardy's inequalities).
Consequently, applying Lemma \ref{CZT} to 
the scalar-valued operator $T = \mathcal{H}_{0,\textrm{loc}}$ we get the following.

\begin{lema} \label{lochil} 
Let $\lambda \in \mathbb{R}$. Then $\mathcal{H}_{\lambda,\emph{loc}}$ is a local Calder\'on-Zygmund
operator, hence, for each $\delta \in \mathbb{R}$, it is bounded on $L^p(\mathbb{R}_+,x^\delta dx)$,
$1<p<\infty$, and from $L^1(\mathbb{R}_+,x^\delta dx)$ into $L^{1,\infty}(\mathbb{R}_+,x^\delta dx)$. 
\end{lema}

We are now prepared to give the main proof.

\subsection*{Proof of Theorem \ref{triesz}} 
We shall justify the sufficiency parts first. To this end assume that $f\in
C_c^{\infty}(\mathbb{R}_+)$. Similarly as in the proof of Theorem \ref{tmax},
we split the Riesz operator
\begin{align} \label{TR1}
R_\lambda f(x) & =
\bigg\{\int_0^{x/2}+ \,\pv \int_{x/2}^{2x}+\int_{2x}^\infty\bigg\}R_\lambda(x,y)f(y)\,
    d\mu_{\lambda}(y)\\
& \equiv R_{\lambda,1}f(x)+R_{\lambda,2}f(x)+R_{\lambda,3}f(x).\nonumber
\end{align}
The operators $R_{\lambda,1}$, $R_{\lambda,2}$ and $R_{\lambda,3}$ 
will be analyzed separately. Note that only the diagonal part $R_{\lambda,2}$
is given by a singular integral.

Applying Lemma \ref{leriesz} we obtain, for all $x>0$,
\begin{align*}
|R_{\lambda,1}f(x)| & \le C_{\lambda} \frac{1}{x^{2\lambda+1}}\int_0^{x/2}|f(y)|y^{2\lambda}dy
\le C_{\lambda} H_0^{2\lambda}|f|(x)\\
|R_{\lambda,3}f(x)|& \le C_{\lambda} x\int_{2x}^\infty {|f(y)|}\frac{dy}{y^2}
\le C_{\lambda} H_{\infty}^1 |f|(x),
\end{align*}
and also 
$$
\big|R_{\lambda,2}f(x)- \mathcal{H}_{\lambda, \textrm{loc}}f(x)
\big| \le C_{\lambda} \int_{x/2}^{2x}\frac{1}{y}\Big(1+\log\frac{xy}{(x-y)^2}\Big)|f(y)|\,dy.
$$
Note that the operator $N$ defined by
$$
Nf(x)=\int_{x/2}^{2x}\frac{1}{y}\Big(1+\log\frac{xy}{(x-y)^2}\Big)f(y)\,dy, \qquad x>0,
$$
and occurring above, is bounded on $L^p(\mathbb{R}_+,x^{\delta}dx)$ for
each $\delta \in \mathbb{R}$ and each $1\le p < \infty$. Indeed, observe first that the
integral defining $N \boldsymbol{1}(x)$
is finite and in fact does not depend on $x>0$. Then, using Jensen's inequality and
changing the order of integration, we get
\begin{align*}
\int_0^\infty|Nf(x)|^p x^{\delta}dx &\le C_p \int_0^\infty
x^\delta\int_{x/2}^{2x}\frac{1}{y}\Big(1+\log\frac{xy}{(x-y)^2}\Big)|f(y)|^p\, dy dx\\
&\le C_p \int_0^\infty
|f(y)|^p\int_{y/2}^{2y} \frac{1}{x}\Big(1+\log\frac{xy}{(x-y)^2}\Big)x^{\delta}{dx} dy\\
&\le C_{p,\delta} \int_0^\infty|f(y)|^py^\delta dy.
\end{align*}

Taking into account Lemmas \ref{Hardy0} and \ref{Hardyinfty}, the above facts and 
Lemma \ref{lochil}, we conclude the following mapping properties of ${R}_{\lambda}$ 
considered on the space $(\mathbb{R}_+,x^\delta dx)$.
For $1<p<\infty$, $R_{\lambda}$ is of strong type $(p,p)$ if $-1-p<\delta<(2\lambda+1)p-1$.
If $-2\le \delta \le 2\lambda$ then $R_{\lambda}$ is of weak type $(1,1)$. Finally,
$R_{\lambda}$ is of restricted weak type $(p,p)$ when $1<p<\infty$ and
$-p-1 \le \delta \le (2\lambda+1)p-1$. 
These properties combined together justify the sufficiency parts of Theorem \ref{triesz}.

The necessity parts of Theorem \ref{triesz} will be justified once we show the following
statements (we assume that $1\le p<\infty$ and the underlying space is $(\mathbb{R}_+,x^\delta dx)$).
\begin{itemize}
\item[(A)]
If $R_{\lambda}$ is of restricted weak type $(p,p)$ then $-p-1\le \delta \le (2\lambda+1)p-1$.
\item[(B1)]
$R_{\lambda}$ is not of weak type $(p,p)$ when $p>1$ and $\delta = -p-1$.
\item[(B2)]
$R_{\lambda}$ is not of weak type $(p,p)$ when $p>1$ and $\delta = (2\lambda+1)p-1$.
\item[(C)]
$R_{\lambda}$ is not of strong type $(1,1)$ for $-2 \le \delta \le 2\lambda$.
\end{itemize}
Item (A) can be concluded immediately by a standard interpolation argument,
once we prove the other items. For item (B1) observe that, 
in view of the above considerations and \eqref{TR1},
$R_{\lambda}$ is of weak type $(p,p)$ for $\delta=-p-1$, $p>1$, if and only if
$R_{\lambda,3}$ has the same property. Moreover, by Lemma \ref{l2riesz}, if
$R_{\lambda,3}$ has this property then the operator $H^1_{\infty}$ also has it.
But it is known (see \cite[Theorem 5]{AM}) that $H^1_\infty$ fails to be
of weak type $(p,p)$ for $\delta = -p-1$, $p>1$. Therefore the same negative result
holds for $R_{\lambda}$; this gives (B1). Treatment of (B2) is similar:
$R_{\lambda}$ is of weak type $(p,p)$ for $\delta = (2\lambda+1)p-1$, $p>1$,
if and only if $R_{\lambda,1}$ has this property. Then, using Lemma \ref{l2riesz},
we infer that the property for $R_{\lambda,1}$ implies the same for $H_0^{2\lambda}$.
However, it is known (see \cite[Theorem 1]{Anker}) that $H_0^{2\lambda}$ does not possess
the property in question, thus (B2) follows.

It remains to show (C). For $0<\varepsilon<1/4$ consider the function
$f_\varepsilon(x)=x^{-\lambda}\chi_{(1,1+\varepsilon)}(x)$, $x>0$.
By the estimate in the diagonal region from Lemma \ref{leriesz} we have
\begin{align*}
\int_0^\infty |R_\lambda f_\varepsilon(x)|x^\delta dx
& \ge \int_{1+2\varepsilon}^2
|R_\lambda f_\varepsilon(x)|x^\delta dx\\
& = \int_{1+2\varepsilon}^2\Big|\int_{x/2}^{2x} R_\lambda(x,y) f_\varepsilon (y)\,d\mu_{\lambda}(y)
\Big|x^\delta dx\\
& \ge \frac{1}{\pi}
\int_{1+2\varepsilon}^2x^{\delta-\lambda}\Big|\int_1^{1+\varepsilon}\frac{1}{y-x}
dy\Big|\,dx-C_{\lambda} \int_{1+2\varepsilon}^{2}
Nf_{\varepsilon}(x) \, x^\delta dx.
\end{align*}
Then, if $R_\lambda$ were bounded on $L^1(\mathbb{R}_+,x^\delta dx)$ we would have
\begin{align*}
\int_{1+2\varepsilon}^2x^{\delta-\lambda}\Big|\int_1^{1+\varepsilon}\frac{1}{y-x}dy\Big|dx
&\le C_{\lambda} \big(\|R_{\lambda}f_\varepsilon\|_{L^1(\mathbb{R}_+,x^\delta dx)}
+\|Nf_\varepsilon\|_{L^1(\mathbb{R}_+,x^\delta dx)}\big)\\
& \le C_{\lambda,\delta} \|f_\varepsilon\|_{L^1(\mathbb{R}_+,x^\delta dx)}
\end{align*}
and hence it would follow that
\begin{align*}
\varepsilon \ge c_{\lambda,\delta} \|f_\varepsilon\|_{L^1(\mathbb{R}_+,x^\delta dx)}
\ge c_{\lambda,\delta}\int_{1+2\varepsilon}^2\log\frac{x-1}{x-1-\varepsilon}dx
= c_{\lambda,\delta}\Big(\varepsilon 
\log\frac{1}{4\varepsilon}+(1-\varepsilon)\log \frac{1}{1-\varepsilon}\Big).
\end{align*}
But the inequality between the outer expressions cannot hold with $c_{\delta,\lambda}$ independent
of $\varepsilon\in (0,1/4)$,  as can be seen immediately by letting $\varepsilon \to 0^+$.
Thus $R_\lambda$ is not bounded on $L^1(\mathbb{R}_+,x^\delta dx)$ and (C) is justified.

The proof of Theorem \ref{triesz} is now complete.
\qed

\begin{remark} 
The boundedness properties of the operator $\widetilde{R}_\lambda$ related to
$\widetilde{\Delta}_\lambda$ can be proved by a careful analysis of the
proof of Theorem \ref{triesz}. Split $\widetilde{R}_\lambda$ into three parts,
$$
\widetilde{R}_\lambda f(x)=\bigg\{\int_0^{x/2} + \pv \int_{x\slash 2}^{2x} + \int_{2x}^{\infty} \bigg\} 
\widetilde{R}_\lambda(x,y)f(y)\,dy
\equiv \widetilde{R}_{\lambda,1}f(x)+\widetilde{R}_{\lambda,2}f(x)+\widetilde{R}_{\lambda,3}f(x).
$$
Since $\widetilde{R}_\lambda(x,y)=(xy)^\lambda
R_\lambda(x,y)$, with the aid of Lemmas \ref{leriesz} and \ref{l2riesz} the
operators $\widetilde{R}_{\lambda,1}$ and
$\widetilde{R}_{\lambda,3}$ can be controlled above and below by
$H_0^\lambda$ and $H_\infty^{\lambda+1}$, respectively. Moreover,
$\widetilde{R}_{\lambda,2}$ is a local Calder\'on-Zygmund operator
with respect to Lebesgue measure; this follows from the
corresponding property for $R_{\lambda,2}$ (see Lemmas \ref{leriesz} and \ref{lochil}). 
Thus we obtain similar results for $\widetilde{R}_\lambda$ to those stated in
Theorem \ref{triesz}, just replacing the role of the interval
$(-p-1,(2\lambda+1)p-1)$ by $(-(\lambda+1)p-1,(\lambda+1)p-1)$.
\end{remark}

\section{The heat integral square function} \label{sec:g}

This section is devoted to the proof of Theorem \ref{tg1}.
Recall that the square function we take into account is given by
$$
g_{\lambda}(f)(x) = \bigg(\int_0^\infty t\Big|\frac{\partial}{\partial t}\int_0^\infty
W_t^\lambda(x,y)f(y)\, d\mu_{\lambda}(y)\Big|^2dt\bigg)^{1/2}, \qquad x>0.
$$
We will need several technical results, one of them being the following important estimate.

\begin{lema} \label{leg1} Let $\lambda>-1/2$. There exists a constant $C_\lambda$ such that
for all $x,y>0$
$$
\bigg(\int_0^\infty t \bigg| \frac{\partial}{\partial t} W_t^\lambda(x,y)
-\chi_{\{0<x/2<y<2x\}}(xy)^{-\lambda}\frac{\partial}{\partial t}
\mathcal{W}_t(x,y)\bigg|^2 dt\bigg)^{1/2}\le C_{\lambda} \big(\max\{x,y\}\big)^{-2\lambda-1}.
$$
\end{lema}

\begin{proof}
We first show that the estimate holds in the diagonal region $0< x\slash 2 \le y \le 2x$.
We shall consider two cases determined by the asymptotics at $0^+$ and $\infty$ of the Bessel
function involved.\\
\noindent \textit{Case 1}: $xy \ge t$. Observe that
\begin{align*}
\frac{\partial}{\partial t}W_t^\lambda(x,y) & =
\frac{\partial}{\partial t}\Big(\frac{1}{\sqrt{2t}}e^{-(x-y)^2/4t}\Big) (xy)^{-\lambda}e^{-xy/2t}
\Big(\frac{xy}{2t}\Big)^{1/2} I_{\lambda-1/2}\Big(\frac{xy}{2t}\Big)\\
&\quad + (xy)^{-\lambda} \frac{\partial}{\partial t}\bigg(e^{-xy/2t}
\Big(\frac{xy}{2t}\Big)^{1/2} I_{\lambda-1/2}\Big(\frac{xy}{2t}\Big)\bigg)
\frac{1}{\sqrt{2t}}e^{-(x-y)^2/4t}\\
& \equiv E_{\lambda,1}(t,x,y)+E_{\lambda,2}(t,x,y).
\end{align*}
Then, according to \eqref{P4}, we get
\begin{align*}
E_{\lambda,1}(t,x,y)&= \frac{(xy)^{-\lambda}}{\sqrt{2\pi}}\frac{\partial}{\partial
t}\Big(\frac{1}{\sqrt{2t}}e^{-(x-y)^2/4t}\Big)\bigg(1+\mathcal{O}\Big(\frac{t}{xy}\Big)\bigg)\\
&= {(xy)^{-\lambda}} \frac{\partial}{\partial
t}\mathcal{W}_t(x,y)-(xy)^{-\lambda}\mathcal{O}\Big(\frac{t}{xy}\Big)\Big(
\frac{1}{(2t)^{3/2}} - \frac{1}{(2t)^{5/2}} (x-y)^2 \Big) e^{-(x-y)^2/4t} \\
& \equiv {(xy)^{-\lambda}} \frac{\partial}{\partial t} \mathcal{W}_t(x,y) +
E_{\lambda,1}^1(t,x,y)+E_{\lambda,1}^2(t,x,y).
\end{align*}
Integrating the two last terms in $t \le xy$ we obtain
\begin{align*}
\bigg(\int_0^{xy}t|E_{\lambda,1}^1(t,x,y)|^2 dt\bigg)^{1/2}
& \le C_{\lambda}(xy)^{-\lambda-1}\bigg(\int_0^{xy}e^{-(x-y)^2/2t}dt\bigg)^{1/2}\\
&\le  C_{\lambda} (xy)^{-\lambda-1/2}, 
\end{align*}
\begin{align*}
\bigg(\int_0^{xy}t|E_{\lambda,1}^2(t,x,y)|^2dt\bigg)^{1/2}
& \le C_{\lambda} (xy)^{-\lambda-1}\bigg(\int_0^{xy}\bigg(\frac{(x-y)^2}{t}\bigg)^2
 e^{-(x-y)^2/2t}dt\bigg)^{1/2}\\
&\le C_{\lambda} (xy)^{-\lambda-1/2}. 
\end{align*}
Next, we analyze $E_{\lambda,2}(t,x,y)$. By using (\ref{P2}) and then (\ref{P4}) it follows that
\begin{align*}
E_{\lambda,2}(t,x,y) & =(xy)^{-\lambda} \frac{\partial}{\partial t} \bigg(e^{-xy/2t}
\Big(\frac{xy}{2t}\Big)^{1/2} I_{\lambda-1/2}\Big(\frac{xy}{2t}\Big)\bigg)
\frac{1}{\sqrt{2t}}e^{-(x-y)^2/4t}\\
&=\frac{(xy)^{-\lambda+1}}{2\sqrt{2}t^{5/2}}e^{-(x-y)^2/4t}
\bigg(e^{-xy/2t}\Big(\frac{xy}{2t}\Big)^{1/2}I_{\lambda-1/2}\Big(\frac{xy}{2t}\Big)\\
&\quad -e^{-xy/2t}\frac{2\lambda t}{xy}\Big(\frac{xy}{2t}\Big)^{1/2}
I_{\lambda-1/2}\Big(\frac{xy}{2t}\Big)
- e^{-xy/2t}\Big(\frac{xy}{2t}\Big)^{1/2} I_{\lambda+1/2}\Big(\frac{xy}{2t}\Big)\bigg)\\
&= \frac{(xy)^{-\lambda+1}}{4\sqrt{\pi} t^{5/2}}e^{-(x-y)^2/4t}
\bigg\{\bigg(1-c(\lambda)\frac{t}{xy}+ \mathcal{O}\Big(\frac{t^2}{(xy)^2}\Big)\bigg)\\
&\quad -\frac{2\lambda t}{xy}\bigg(1-c(\lambda)\frac{t}{xy}+\mathcal{O}\Big(\frac{t^2}{(xy)^2}\Big)\bigg)
-\bigg(1- c(\lambda+1)\frac{t}{xy}+\mathcal{O}\Big(\frac{t^2}{(xy)^2}\Big)\bigg)\bigg\},
\end{align*}
where $c(\nu)=[\nu-1/2,1]=(\nu-1/2)^2-1/4$. Hence, due to the
occurring cancellations,
$$
|E_{\lambda,2}(t,x,y)|=
\frac{(xy)^{-\lambda+1}}{t^{5/2}}e^{-(x-y)^2/4t}\mathcal{O}\Big(\frac{t^2}{(xy)^2}\Big).
$$
Consequently,
$$
\bigg(\int_0^{xy}t|E_{\lambda,2}(t,x,y)|^2dt\bigg)^{1/2} \le
C_{\lambda}(xy)^{-\lambda-1}\bigg(\int_0^{xy}e^{-(x-y)^2/2t}dt\bigg)^{1/2}
\le C_{\lambda} (xy)^{-\lambda-1/2}. 
$$
In view of this and the previous estimates we conclude that
$$
\bigg(\int_0^{xy} t \Big|\frac{\partial}{\partial t}W_t^\lambda(x,y)
-(xy)^{-\lambda}\frac{\partial}{\partial t}\mathcal{W}_t(x,y)\Big|^2 dt\bigg)^{1/2}
\le C_{\lambda} (xy)^{-\lambda-1/2}, \qquad x,y>0, 
$$
which finishes Case 1.\\
\noindent \textit{Case 2}: $xy < t$. Observe that (\ref{P2}) leads to
\begin{align} \nonumber
\frac{\partial}{\partial t} W_t^\lambda(x,y) & = 2^{-\lambda-1/2} e^{-(x^2+y^2)/4t}
\bigg(-\frac{\lambda+1/2}{t^{\lambda+3/2}}\Big(\frac{xy}{2t}\Big)^{-\lambda+1/2}
I_{\lambda-1/2}\Big(\frac{xy}{2t}\Big)\\
&\quad -\frac{xy}{2t^{\lambda+5/2}}\Big(\frac{xy}{2t}\Big)^{-\lambda+1/2}
I_{\lambda+1/2}\Big(\frac{xy}{2t}\Big) 
+\frac{x^2+y^2}{4t^{\lambda+5/2}}\Big(\frac{xy}{2t}\Big)^{-\lambda+1/2}
I_{\lambda-1/2}\Big(\frac{xy}{2t}\Big)\bigg).\label{41}
\end{align}
Then \eqref{P3} implies
\begin{equation}\label{42}
\Big|\frac{\partial}{\partial t}W_t^\lambda(x,y)\Big| \le
\frac{C_{\lambda}}{t^{\lambda+3/2}} \Big(1+
\frac{(xy)^2}{t^2}+\frac{x^2+y^2}{t}\Big)e^{-(x^2+y^2)/4t} \le
\frac{C_{\lambda}}{t^{\lambda+3/2}}e^{-(x^2+y^2)/8t},
\end{equation}
hence we can write
\begin{align} \nonumber
\bigg(\int_{xy}^\infty t\Big|\frac{\partial}{\partial t}W_t^\lambda(x,y)\Big|^2 dt\bigg)^{1/2}
&\le C_{\lambda} \bigg(\int_{xy}^\infty \frac{1}{t^{2\lambda+2}}e^{-(x^2+y^2)/4t}dt\bigg)^{1/2}\\
\label{77} &\le C_{\lambda}
\frac{1}{(x^2+y^2)^{\lambda+1/2}}\bigg(\int_0^{(x^2+y^2)/4xy}u^{2\lambda}e^{-u}du\bigg)^{1/2}\\
&\le C_{\lambda} (xy)^{-\lambda-1/2}. \nonumber
\end{align}
Moreover, by a straightforward computation,
$$
\Big|\frac{\partial}{\partial t}\mathcal{W}_t(x,y)\Big| \le
\frac{C}{t^{3/2}}e^{-(x-y)^2/4t} \Big(1+\frac{(x-y)^2}{t}\Big) \le \frac{C}{t^{3/2}}e^{-(x-y)^2/8t} 
$$
and therefore
\begin{align} \nonumber
\bigg(\int_{xy}^\infty
t\Big|(xy)^{-\lambda}\frac{\partial}{\partial t}\mathcal{W}_t(x,y)\Big|^2 dt \bigg)^{1\slash 2}
& \le C (xy)^{-\lambda} \bigg(\int_{xy}^\infty \frac{1}{t^2}e^{-(x-y)^2/4t}dt\bigg)^{1/2} \\
\label{410} & \le C(xy)^{-\lambda-1/2}. 
\end{align}

Now, combining the estimates of Cases 1 and 2, we conclude that for $0< x/2 \le y \le 2x$,
$$
\bigg(\int_0^\infty t\Big|\frac{\partial}{\partial t}
W_t^\lambda(x,y)- (xy)^{-\lambda}\frac{\partial}{\partial t}\mathcal{W}_t(x,y)\Big|^2 dt\bigg)^{1/2}
\le C_{\lambda} (xy)^{-\lambda-1\slash 2} \le C_{\lambda} \big(\max\{x,y\}\big)^{-2\lambda-1}.
$$

It remains to prove the relevant bound in the off-diagonal region. By symmetry of the
kernel $W_t^{\lambda}(x,y)$ it is enough to focus on the cone $0 < y < x\slash 2$.
Similarly as above, we consider two cases.\\
\noindent \textit{Case 1}: $xy \ge t$.
According to \eqref{41} and \eqref{P4},
\begin{align*}
\frac{\partial}{\partial t} W_t^\lambda(x,y)
& = \frac{2^{-\lambda-1/2}}{\sqrt{2\pi}}e^{-(x-y)^2/4t}\bigg\{-\frac{\lambda+1/2}{2^{-\lambda}t^{3/2}}
(xy)^{-\lambda}\bigg(1+\mathcal{O}\Big(\frac{t}{xy}\Big)\bigg)\\
&\quad -\frac{(xy)^{-\lambda+1}}{2^{-\lambda+1}t^{5/2}}
\bigg(1+\mathcal{O}\Big(\frac{t}{xy}\Big)\bigg)
+\frac{(xy)^{-\lambda}(x^2+y^2)}{2^{-\lambda+2}t^{5/2}}
\bigg(1+\mathcal{O}\Big(\frac{t}{xy}\Big)\bigg)\bigg\}.
\end{align*}
This implies
\begin{align*}
\Big|\frac{\partial}{\partial t}W_t^\lambda(x,y)\Big| & \le
C_{\lambda} \frac{(xy)^{-\lambda}}{t^{1/2}}e^{-x^2/16t}\Big(\frac{1}{t}+\frac{1}{xy}+
\frac{xy}{t^2} + \frac{x^2+y^2}{t^2}+\frac{x^2+y^2}{xyt}\Big) \\
& \le C_{\lambda} \frac{(xy)^{-\lambda}x^2}{t^{5/2}}e^{-x^2/16t},
\end{align*}
since by the present assumptions on $x,y$ and $t$ we have
$$
\frac{1}{xy} \le \frac{1}{t} \le \frac{x^2+y^2}{txy} \le \frac{x^2+y^2}{t^2} \le 2\frac{x^2}{t^2}.
$$
Using again the assumptions $xy \ge t$ and $y< x\slash 2$ it follows that
\begin{align*}
\bigg(\int_0^{xy}t\Big|\frac{\partial}{\partial t}W_t^\lambda(x,y)\Big|^2 dt\bigg)^{1/2}
& \le C_{\lambda}
\bigg(\int_0^{xy}t\Big|(xy)^{-\lambda}\frac{x^2}{t^{5/2}}e^{-x^2/16t}\Big|^2 dt\bigg)^{1/2} \\
& \le C_{\lambda} \frac{1}{x^{2\lambda+3}} \bigg(xy\int_{0}^{xy}
\Big(\frac{x^2}{t}\Big)^{2\lambda+5}e^{-x^2/8t} dt\bigg)^{1/2}\\
& \le C_{\lambda} \frac{xy}{x^{2\lambda+3}},
\end{align*}
and the last quantity can be easily estimated from above by $C_{\lambda} x^{-2\lambda-1}$.

\noindent \textit{Case 2}: $xy < t$.
Notice that by \eqref{77} we have
$$
\bigg(\int_{xy}^\infty t\Big|\frac{\partial}{\partial t}W_t^\lambda(x,y)\Big|^2 dt\bigg)^{1/2}
\le \frac{C_{\lambda}}{(x^2+y^2)^{\lambda+1/2}} \le C_{\lambda} x^{-2\lambda-1}.
$$

Putting the two cases together produces
$$
\bigg(\int_0^\infty t\Big|\frac{\partial}{\partial t}W_t^\lambda(x,y)\Big|^2 dt\bigg)^{1/2}
\le C_{\lambda} x^{-2\lambda-1}, \qquad 0<y<x/2,
$$
which is precisely what we needed. 

The proof of Lemma \ref{leg1} is finished.
\end{proof}

\begin{lema} \label{4.2}
Let $\lambda>-1/2$. There exist (positive) constants $a=a_{\lambda}$ and $c = c_{\lambda}$ such that
$$
\frac{\partial}{\partial t}W_t^\lambda(x,y) \le - \frac{c}{t^{\lambda+3/2}}
$$
if either $0 < x,y < a$ and $t\ge 1$ or $0< y < x$ and $x^2\slash t \le a$.
\end{lema}

\begin{proof}
The conclusion follows in a straightforward manner by combining \eqref{41} and \eqref{P3}.
\end{proof}

Define the auxiliary square function
\begin{equation} \label{gvlocdf}
\mathfrak{g}_{\lambda,\textrm{loc}}(f)(x)=\bigg(\int_0^\infty
t\Big|\int_{x/2}^{2x}(xy)^{-\lambda}\frac{\partial}{\partial
t}\mathcal{W}_t(x,y)f(y)\,d\mu_{\lambda}(y)\Big|^2 dt\bigg)^{1/2}, \qquad x>0,
\end{equation}
which is the local part of a modification of the classic vertical square function.

\begin{lema} \label{gVloc}
Let $\lambda \in \mathbb{R}$. Then, for each $\delta\in \mathbb{R}$,
$\mathfrak{g}_{\lambda,\textrm{loc}}$ is bounded on
$L^p(\mathbb{R}_+,x^\delta dx)$, $1< p< \infty$,
and from $L^1(\mathbb{R}_+,x^\delta dx)$ into $L^{1,\infty}(\mathbb{R}_+,x^\delta dx)$.
\end{lema}

\begin{proof}
Observe that the mapping $\mathfrak{g}_{0,\textrm{loc}}$
is bounded on $L^2(\mathbb{R}_+,dx)$. This can be easily verified
by invoking the well-known fact that the classical vertical $g$-function $\mathfrak{g}$
(restricted to $\mathbb{R}_+$ and given by a formula analogous to that defining
$\mathfrak{g}_{0,\textrm{loc}}$, but with the integration in $y$ from $0$ to $\infty$) 
is bounded on $L^2(\mathbb{R}_+,dx)$ and then using classic Hardy's inequalities.

From the $L^2$-boundedness of $\mathfrak{g}_{0,\textrm{loc}}$ we infer that the
operator $\mathcal{G}$ assigning to an $f \in L^2(\mathbb{R}_+,dx)$ the function
$$
\mathbb{R}_+ \ni x \mapsto \mathcal{G}f(x) = \bigg\{ \int_{x\slash 2}^{2x}
     \frac{\partial}{\partial t} \mathcal{W}_t(x,y) f(y)\, dy\bigg\}_{t>0},
$$
is bounded from $L^2(\mathbb{R}_+,dx)$ to the Bochner-Lebesgue space
$L^2_{L^2(tdt)}(\mathbb{R}_+,dx)$. Moreover, $\mathcal{G}$ is associated, in the
sense of Lemma \ref{CZT}, with the vector-valued kernel 
$$
K(x,y) = \Big\{\frac{\partial}{\partial t} \mathcal{W}_t(x,y) \Big\}_{t>0}.
$$ 
This follows essentially by the known fact that $\mathfrak{g}$, viewed as a vector-valued
Calder\'on-Zygmund operator, is associated with the same kernel. 
In addition, $K$ satisfies the standard estimates from Lemma \ref{CZT}
(even for all $x,y\in \mathbb{R}$), as is known and not hard to check.

Taking into account the above facts and applying Lemma \ref{CZT} with $T = \mathcal{G}$
we see that the mapping $f(x) \mapsto x^{-\lambda}\mathcal{G}(y^{\lambda}f)(x)$ 
is a local vector-valued Calder\'on-Zygmund operator, hence, 
given arbitrary $\delta \in \mathbb{R}$, it extends to a bounded operator from
$L^p(\mathbb{R}_+,x^{\delta}dx)$ to $L^p_{L^2(tdt)}(\mathbb{R}_+,x^{\delta}dx)$ for any
$1<p<\infty$, and to a bounded operator from $L^1(\mathbb{R}_+,x^{\delta}dx)$ to
$L^{1,\infty}_{L^2(tdt)}(\mathbb{R}_+,x^{\delta}dx)$.

Finally, notice that these boundedness results imply precisely the desired mapping properties of
$\mathfrak{g}_{\lambda,\textrm{loc}}$ (the fact that the vector-valued bounded extensions correspond 
to $\mathfrak{g}_{\lambda,\textrm{loc}}$ given by \eqref{gvlocdf} follows by a standard density argument). 
\end{proof}

Having the above results we are ready to prove Theorem \ref{tg1}. 

\subsection*{Proof of Theorem \ref{tg1}} 
Using the triangle inequality for the norm $\|\cdot\|_{L^2(tdt)}$ and then applying Minkowski's
integral inequality we get
\begin{align*}
g_{\lambda}(f)(x)
& \le
\bigg\{\int_0^{x/2} + \int_{2x}^{\infty}\bigg\} \bigg(\int_0^\infty t\Big|\frac{\partial}{\partial t}
W_t^\lambda(x,y)\Big|^2dt\bigg)^{1/2} |f(y)| d\mu_{\lambda}(y)\\
& \quad + \int_{x/2}^{2x} \bigg(\int_0^\infty t\Big|\frac{\partial}{\partial
t}W_t^\lambda(x,y)-{(xy)^{-\lambda}}\frac{\partial}{\partial t}\mathcal{W}_t(x,y)\Big|^2dt\bigg)^{1/2}
|f(y)| d\mu_{\lambda}(y)\\ &\quad + \bigg(\int_0^\infty
t\Big|\int_{x/2}^{2x}(xy)^{-\lambda}\frac{\partial}{\partial
t}\mathcal{W}_t(x,y)f(y) d\mu_{\lambda}(y)\Big|^2dt\bigg)^{1/2}.
\end{align*}
Then by Lemma \ref{leg1} it follows that
$$
g_{\lambda}(f)(x) \le C_{\lambda}\big(H_0^{2\lambda}|f|(x)+H_\infty^0|f|(x)+\mathcal{N}|f|(x)
    +\mathfrak{g}_{\lambda,\textrm{loc}}(f)(x)\big), \qquad x>0,
$$
where $\mathcal{N}$ denotes the operator
$$
\mathcal{N}f(x)=\int_{x/2}^{2x}\frac{f(y)}{y}dy, \qquad x>0.
$$
Note that $\mathcal{N}$ is bounded on $L^p(\mathbb{R}_+,x^\delta
dx)$ for each $1\le p \le \infty$ and each $\delta \in \mathbb{R}$
(this can be easily verified by means of Jensen's inequality and
Fubini's theorem, as in the case of the operator $N$ emerging in the
proof of Theorem \ref{triesz}). Furthermore, by Lemma \ref{gVloc},
$\mathfrak{g}_{\lambda,\textrm{loc}}$ is bounded on $L^p(\mathbb{R}_+,x^\delta dx)$, $1< p< \infty$,
and from $L^1(\mathbb{R}_+,x^\delta dx)$ into $L^{1,\infty}(\mathbb{R}_+,x^\delta dx)$, for
each $\delta\in \mathbb{R}$.

Then, similarly as in the proof of Theorem \ref{triesz}, by means of
Lemmas \ref{Hardy0} and \ref{Hardyinfty} we conclude that
$g_{\lambda}$ is bounded on $L^p(\mathbb{R}_+,x^\delta dx)$ when $1<
p< \infty$ and $-1<\delta<(2\lambda+1)p-1$, it is bounded from
$L^1(\mathbb{R}_+,x^\delta dx)$ into
$L^{1,\infty}(\mathbb{R}_+,x^\delta dx)$ when $-1<\delta\le
2\lambda$, and finally it is of restricted weak type $(p,p)$ with
respect to $(\mathbb{R}_+,x^\delta dx)$ when $1\le p<\infty$ and
$-1<\delta\le (2\lambda+1)p-1$. These facts justify the sufficiency parts in Theorem \ref{tg1}.

We now prove the necessity parts in Theorem \ref{tg1}. This task will be done once we show the 
following statements (we assume that $\lambda>-1\slash 2$, $1\le p < \infty$ and the underlying 
space is $(\mathbb{R}_+,x^{\delta}dx)$).
\begin{itemize}
\item[(A)]
If $g_{\lambda}$ is of restricted weak type $(p,p)$ then $-1<\delta\le (2\lambda+1)p-1$.
\item[(B)]
$g_{\lambda}$ is not of weak type $(p,p)$ when $p>1$ and $\delta = (2\lambda+1)p-1$.
\item[(C)]
$g_{\lambda}$ is not of strong type $(1,1)$ if $-1 < \delta \le 2\lambda$.
\end{itemize}

Let $f=\chi_{(a\slash 2,a)}$, where $a$ is the constant from Lemma \ref{4.2}. Then, according to that lemma,
\begin{align*}
g_{\lambda}(f)(x) 
&\ge \bigg(\int_1^\infty t\Big|\int_{a/2}^a\frac{\partial}{\partial
t}W_t^\lambda(x,y)\, d\mu_{\lambda}(y)\Big|^2dt\bigg)^{1/2} \\
& \ge c_{\lambda} \bigg(\int_1^\infty
\frac{1}{t^{2\lambda+2}}\Big|\int_{a\slash 2}^ay^{2\lambda}dy\Big|^2dt\bigg)^{1/2} \\
& = c_{\lambda}, \qquad 0 < x < a.
\end{align*}
Suppose that $g_\lambda$ can be extended from $L^2(\mathbb{R}_+,d\mu_{\lambda})$ to 
a restricted weak type $(p,p)$ operator on $(\mathbb{R}_+,x^\delta dx)$.
For sufficiently small $\gamma>0$ we then have 
$$
\int_0^a x^\delta dx \le C_{p,\lambda} \gamma^{-p} {\|f\|^p_{L^p(\mathbb{R}_+,x^\delta dx)}}
    \le C_{p,\lambda} \gamma^{-p},
$$
which implies $\delta>-1$.

On the other hand, using again Lemma \ref{4.2}, we get
\begin{align} \nonumber
g_\lambda(f)(x) & \ge c_\lambda\bigg(\int_{x^2/a}^\infty
\frac{1}{t^{2\lambda+2}}\Big|\int_0^1 f(y) \, d\mu_{\lambda}(y)\Big|^2dt \bigg)^{1/2} \\
& \ge c_\lambda x^{-2\lambda-1}\int_0^1 f(y)\, d\mu_{\lambda}(y), \qquad x>1, \label{we2}
\end{align}
for any nonnegative measurable function $f$ on $\mathbb{R}_+$. Thus
taking $f = \chi_{(1\slash 2,1)}$ we have $g_{\lambda}(f)(x) \ge
c_{\lambda}x^{-2\lambda-1}$, provided that $x>1$. Consequently, if
$g_\lambda$ can be extended from $L^2(\mathbb{R}_+,d\mu_{\lambda})$
to a restricted weak type $(p,p)$ operator on $(\mathbb{R}_+,x^\delta dx)$, then
$$
\int_1^{\gamma^{-1/(2\lambda+1)}}x^\delta dx \le
C_{p,\lambda} \gamma^{-p} \|f\|^p_{L^p(\mathbb{R}_+,x^\delta dx)} \le C_{p,\lambda} \gamma^{-p}
$$
for $\gamma>0$ small enough. This implies $\delta \le (2\lambda+1)p-1$ and (A) is justified.

In order to show (B) we use the estimate \eqref{we2} and proceed as in the proof of
the corresponding result in Theorem \ref{tmax}.

Considering (C), we argue similarly as in the proofs of the parallel properties in
Theorems \ref{tmax} and \ref{triesz}. We shall first see that
$\mathfrak{g}_{\lambda,\textrm{loc}}$ is not bounded on
$L^1(\mathbb{R}_+,x^\delta dx)$. Let $0<\varepsilon<1/2$ and assume
that $1< y< 1+\varepsilon$, $1+2\varepsilon < x < 2$ and $t\ge 2(x-1)^2$. Then we have
$$
\frac{\partial}{\partial t}\mathcal{W}_t(x,y) =
\frac{1}{2\sqrt{\pi}t^{3/2}}\Big(-1+\frac{(x-y)^2}{2t}\Big)e^{-(x-y)^2/4t} \le - c t^{-3/2},
$$
for some $c>0$, because $(x-y)^2/2t\le 9\slash 16$, see \eqref{hasz}. Hence, letting 
$f_\varepsilon = \chi_{(1,1+\varepsilon)}$, it follows that
\begin{align*}
\mathfrak{g}_{\lambda,\textrm{loc}}(f_\varepsilon)(x) & =
\bigg(\int_0^\infty t\Big|\int_{x/2}^{2x}\frac{\partial}{\partial
t}\mathcal{W}_t(x,y)(xy)^{-\lambda} f_\varepsilon(y) \, d\mu_{\lambda}(y) \Big|^2dt\bigg)^{1/2}\\
&\ge c_{\lambda} \bigg(\int_{2(x-1)^2}^\infty\frac{dt}{t^2} \bigg)^{1/2}
\bigg(\int_1^{1+\varepsilon}dy\bigg)\\
&\ge c_{\lambda} \frac{\varepsilon}{x-1}, \qquad x\in (1+2\varepsilon,2).
\end{align*}
Consequently,
$$
\int_0^\infty \mathfrak{g}_{\lambda,\textrm{loc}}(f_\varepsilon)(x) x^\delta dx \ge \int_{1+2\varepsilon}^2
\mathfrak{g}_{\lambda,\textrm{loc}}(f_\varepsilon)(x) x^\delta dx
\ge c_{\lambda} \varepsilon \int_{1+2\varepsilon}^2\frac{x^\delta
dx}{x-1} \ge c_{\lambda,\delta} \varepsilon \log\frac{1}{2\varepsilon}.
$$
Now, if $\mathfrak{g}_{\lambda,\textrm{loc}}$ were bounded on
$L^1(\mathbb{R}_+,x^\delta dx)$ we would have
\begin{equation}\label{412}
c_{\lambda,\delta} \varepsilon \log \frac{1}{2\varepsilon} \le
\int_0^\infty \mathfrak{g}_{\lambda,\textrm{loc}}(f_\varepsilon)(x)
x^\delta dx \le C \int_0^\infty f_\varepsilon(x) x^\delta dx \le
C_{\delta} \varepsilon,
\end{equation}
which obviously cannot hold as $\varepsilon \to 0^+$. 
Thus $\mathfrak{g}_{\lambda,\textrm{loc}}$ is not bounded on
$L^1(\mathbb{R}_+,x^\delta dx)$. On the other hand, by Lemma
\ref{leg1} and Minkowski's integral inequality we see that, for $x>0$,
\begin{align*}
\mathfrak{g}_{\lambda,\textrm{loc}}(f_\varepsilon)(x)
& \le \bigg(\int_0^\infty t\bigg|\int_{x/2}^{2x}
\bigg(\frac{\partial}{\partial t} \Big((xy)^{-\lambda}
\mathcal{W}_t(x,y)\Big)-\frac{\partial}{\partial t}
W_t^\lambda(x,y)\bigg)f_\varepsilon(y)\, d\mu_{\lambda}(y)\bigg|^2dt\bigg)^{1/2}\\
&\quad + \bigg(\int_0^\infty
t\bigg|\int_{x/2}^{2x}\frac{\partial}{\partial t}
W_t^\lambda(x,y)f_\varepsilon(y)\,d\mu_{\lambda}(y)\bigg|^2dt\bigg)^{1/2}\\
&\le C_{\lambda}
\bigg\{\int_{x/2}^{2x}\frac{|f_\varepsilon(y)|}{y}dy+\bigg(\int_0^\infty
t\bigg|\int_{x/2}^{2x}\frac{\partial}{\partial t}
W_t^\lambda(x,y)f_\varepsilon(y)\, d\mu_{\lambda}(y)\bigg|^2dt\bigg)^{1/2}\bigg\}.
\end{align*}
This gives 
$$
\mathfrak{g}_{\lambda,\textrm{loc}}(f_\varepsilon)(x) \le
C_{\lambda} \big( \mathcal{N}|f_\varepsilon|(x) +g_{\lambda}(f_\varepsilon)(x) \big), \qquad x \in
(1+2\varepsilon,2).
$$
Hence, taking into account weighted $L^1$-boundedness of the
operator $\mathcal{N}$, if $g_{\lambda}$ were bounded on
$L^1(\mathbb{R}_+,x^\delta dx)$ then \eqref{412} would hold, for all
$0<\varepsilon<1/2$, a contradiction. Thus $g_{\lambda}$ is not
bounded on $L^1(\mathbb{R}_+,x^\delta dx)$ and statement (C) follows.

The proof of Theorem \ref{tg1} is complete.
\qed

\begin{remark} Similar facts to those from Remark \ref{Re1} are true for
the corresponding $g$-function in the $\widetilde{\Delta}_\lambda$-setting.
In particular, the relevant interval for $\widetilde{g}_{\lambda}$ is $(-\lambda p -1, (\lambda + 1)p -1)$.
\end{remark}

\section{Operators related to the Poisson integral}
\label{sec:subord}

Recall that the operators analyzed in the previous sections were defined by means of the heat kernel
$W_t^{\lambda}(x,y)$. In this section we consider the maximal operator and a square function related to the
Poisson kernel $P_t^{\lambda}(x,y)$ associated with $\Delta_{\lambda}$, $\lambda>-1\slash 2$.

First of all, we shall compute $P_t^{\lambda}(x,y)$. By the principle of subordination,
\begin{equation} \label{subord}
P_{t}^{\lambda}(x,y) = \int_{0}^{\infty} W^{\lambda}_{t^2\slash 4u}(x,y) \, \frac{e^{-u}du}{\sqrt{\pi u}}.
\end{equation}
Then an application of the integral formula \eqref{intF} leads to
$$
P_t^{\lambda}(x,y) = \mathcal{C}(\lambda)
    \frac{t}{(x^2+y^2+t^2)^{\lambda+1}} \;{_2F_1}\bigg( \frac{\lambda+1}{2}, \frac{\lambda+2}{2};
        \frac{2\lambda+1}{2}; \Big( \frac{2xy}{x^2+y^2+t^2} \Big)^2 \bigg),
$$
with $\mathcal{C}(\lambda) = {2\pi^{-1\slash 2}\Gamma(\lambda+1)}\slash{\Gamma(\lambda+1\slash 2)}$.
We now transform the above expression in order to see the exact behavior of the Poisson kernel.
Using the formula (cf. \cite[(9.5.3)]{Leb}),
$$
{_2F_1}(\alpha,\beta;\gamma;z) = (1-z)^{\gamma-\alpha-\beta}\;
    {_2F_1}(\gamma-\alpha,\gamma-\beta; \gamma;z),
$$
valid for $z<1$, one gets
$$
P_t^{\lambda}(x,y) = \mathcal{C}(\lambda)
    \frac{t(x^2+y^2+t^2)^{1-\lambda}}{[(x+y)^2+t^2][(x-y)^2+t^2]}
    \;{_2F_1}\bigg( \frac{\lambda}{2}, \frac{\lambda-1}{2};
        \frac{2\lambda+1}{2}; \Big( \frac{2xy}{x^2+y^2+t^2} \Big)^2 \bigg).
$$
We observe that the values of the last ${_2F_1}$ function are separated from $0$ and $\infty$;
this is because the function $z\mapsto {_2F_1}(\lambda\slash 2, (\lambda-1)\slash 2;
\lambda+1\slash 2; z)$ is continuous on $[0,1)$, has value $1$ at $z=0$ (see \cite[(9.1.1)]{Leb}),
its limit as $z\to 1^-$ exists and is positive (cf. \cite[Section 9.3]{Leb}) and, finally,
there are no zeroes in $(0,1)$ since the Poisson kernel is strictly positive (this is of course
a consequence of the same property for the heat kernel). Thus we obtain the following.

\begin{propo} \label{p_ker}
Let $\lambda > -1\slash 2$. There exists a constant $C_{\lambda}$ such that, for all $t,x,y>0$,
$$
\frac{C^{-1}_{\lambda}t}{(x^2+y^2+t^2)^{\lambda}[(x-y)^2+t^2]} \le P_t^{\lambda}(x,y)
    \le \frac{C_{\lambda}t}{(x^2+y^2+t^2)^{\lambda}[(x-y)^2+t^2]}.
$$
\end{propo}

It should not be surprising that in the case $\lambda=0$ the Poisson kernel has a particularly
simple form. Notice that the first parameter in the ${_2F_1}$ expressing $P_{t}^0(x,y)$
is zero, consequently the hypergeometric function is identically equal to $1$ and we get
$$
P_t^0(x,y) = \frac{1}{\pi} \bigg( \frac{t}{(x-y)^2+t^2} + \frac{t}{(x+y)^2+t^2} \bigg).
$$
This shows that $P_t^0(x,y)= \mathcal{P}_t(x,y)+\mathcal{P}_t(x,-y)$, with $\mathcal{P}_t$ being
the classical Poisson kernel. Clearly, the same identity could be concluded immediately from
the analogous, already commented in Section \ref{sec:riesz},
connection between $W^0_t$ and $\mathcal{W}_t$.

Apart from the maximal operator $P_*^{\lambda}$ we consider
the Littlewood-Paley type square function
$$
\mathbb{g}_{\lambda}(f)(x) =\bigg(\int_0^\infty t\Big|\frac{\partial}{\partial t}
P_t^{\lambda}f(x)\Big|^2dt\bigg)^{1/2}.
$$
It turns out that various boundedness results for the two
operators can be concluded, in a straightforward manner, by means of
Theorems \ref{tmax} and \ref{tg1}. A key fact here is
that $P^\lambda_*f$ and $\mathbb{g}_{\lambda}(f)$ can be controlled
pointwise by $W^{\lambda}_*f$ and  $g_{\lambda}(f)$, respectively.

\begin{propo} \label{cppp}
Let $\lambda > -1\slash 2$. Then, for sufficiently regular functions $f$, 
$$
P^{\lambda}_*f(x) \le W^{\lambda}_*f(x), \qquad
    \mathbb{g}_{\lambda}(f)(x) \le \sqrt{2} {g}_{\lambda}(f)(x), \qquad x>0.
$$
\end{propo}

\begin{proof}
By subordination and Fubini's theorem
$$
P^{\lambda}_*f(x) = \sup_{t>0} \bigg| \int_0^{\infty} W^{\lambda}_{t^2\slash 4u} f(x)
    \frac{e^{-u}du}{\sqrt{\pi u}} \bigg| \le \int_0^{\infty} \sup_{t>0}
    \big| W^{\lambda}_{t^2\slash 4u} f(x) \big| \frac{e^{-u}du}{\sqrt{\pi u}} =
    W^{\lambda}_*f(x).
$$
Next, again by subordination, 
\begin{align*}
\mathbb{g}_{\lambda}(f)(x) & = \bigg( \int_0^{\infty} t \bigg|
\frac{\partial}{\partial t} \int_0^{\infty}
    W^{\lambda}_{t^2\slash 4u}f(x) \frac{e^{-u}du}{\sqrt{\pi u}} \bigg|^2 dt \bigg)^{1\slash 2}\\
& = \bigg\| \int_0^{\infty} (\partial_t W^{\lambda})_{t^2\slash 4u}f(x)
     \frac{t}{2u}\frac{e^{-u}du}{\sqrt{\pi u}} \bigg\|_{L^2(tdt)},
\end{align*}
where $(\partial_t W^{\lambda})_{t^2\slash 4u}f(x)$ is the derivative
in $s$ of $W_s^{\lambda}f(x)$, taken at the point $s=t^2\slash 4u$.
Then applying Minkowski's integral inequality and changing the variable we obtain
\begin{align*}
\mathbb{g}_{\lambda}(f)(x) & \le \int_0^{\infty} \big\|
t(\partial_tW^{\lambda})_{t^2\slash 4u} f(x)\big\|_{L^2(tdt)}
    \frac{e^{-u}du}{2u\sqrt{\pi u}} \\
& = \sqrt{2}\int_0^{\infty} \frac{e^{-u}du}{\sqrt{\pi u}}\, \bigg( \int_0^{\infty} s
    \Big| \frac{\partial}{\partial s} W^{\lambda}_s f(x) \Big|^2 ds \bigg)^{1\slash 2}\\
& = \sqrt{2} \, g_{\lambda}(f)(x).
\end{align*}
The proof is complete.
\end{proof}

Thus, by Proposition \ref{cppp} and Theorems \ref{tmax} and
\ref{tg1}, we obtain strong type, weak type and
restricted weak type mapping properties of the Poisson integral
based operators. However, since the subordination principle means
a kind of averaging, it is far from obvious whether these results
are sharp in the sense of admissible powers $\delta$. In the case
of $P^{\lambda}_*$ the precise behavior of the corresponding
kernel is known (Proposition \ref{p_ker}) and is relatively
simple. Hence we can easily obtain the following estimates,
similar to those from Section 3 and now involving the Poisson kernel.

\begin{propo} \label{estPoi} Let $\lambda>-1/2$. There exists $c_{\lambda}>0$ such that
$$
P_t^\lambda(x,y)\ge c_\lambda \frac{t}{t^2+(x-y)^2}, \qquad t\le 1\slash 2, \quad x,y\in (1,2),
$$
and such that, for every nonnegative measurable function $f$ on $(0,\infty)$, we have
$$
P_*^\lambda f(x)\ge c_\lambda\int_0^\infty
\frac{f(y)}{(1+y^2)^{\lambda+1}}\, d\mu_{\lambda}(y),\qquad x\in (0,1),
$$
and
$$
P_*^\lambda f(x)\ge c_\lambda x^{-2\lambda-1}\int_0^x f(y)\, d\mu_{\lambda}(y),\qquad x\in (0,\infty).
$$
\end{propo}

Now, using the arguments from the proofs of Theorem \ref{tmax} and Theorem \ref{tg1} to justify
necessity parts, we obtain sharp results for the Poisson integral maximal operator
stated in Theorem \ref{SPoi}.

Proving similar results for 
$\mathbb{g}_{\lambda}$ requires a deeper and more subtle analysis,
which is beyond the scope of this paper. We only mention that in
order to obtain suitable kernel estimates one has to deal with terms
involving hypergeometric functions with different parameters and
examine essential cancellations occurring between those terms (for
that purpose it seems to be more convenient to use, instead of
${_2F_1}$, the integral representation \eqref{subord} of
$P^{\lambda}_t(x,y)$ in terms of $W^{\lambda}_t(x,y)$).


\section*{Acknowledgments}

The authors gratefully acknowledge the hospitality of Professor
Jos\'e Luis Torrea and the Department of Mathematics at Universidad
Aut\'onoma de Madrid. A part of the investigation of the first, second
and fourth-named authors during the preparation of this paper was made
while they were visiting Universidad Aut\'onoma de Madrid. The
third-named author was partially working on this paper during a
postdoctoral visiting position period at the Department of Mathematics 
of Universidad Aut\'onoma de Madrid, October 2005--March 2006.


\end{document}